\documentclass[11pt]{article}
\usepackage{amsthm,amssymb,amsmath,amscd}
\usepackage[all]{xy}
\usepackage{mathrsfs}
\usepackage{multirow}
\usepackage{color}
\usepackage{mathptm}

\newtheorem{Def}{Definition}[section]
\newtheorem{Prop}[Def]{Proposition}
\newtheorem{Theo}[Def]{Theorem}
\newtheorem{Lem}[Def]{Lemma}
\newtheorem{Koro}[Def]{Corollary}

\newcommand{\E}[4]{{\rm E}_{#1}^{#2}(#3, #4)}
    \newcommand{\Ex}[3]{{\rm E}_{#1}^{#2}(#3)}

\newcommand{\id}{{\rm id}}
\newcommand{\add}{{\rm add}}

\newcommand{\domdim}{{\rm dom.dim}}
\newcommand{\Hom}{{\rm Hom }}

\newcommand{\soc}{{\rm soc}}
\renewcommand{\top}{{\rm top}}
\newcommand{\StHom}{{\rm \underline{Hom} }}

\newcommand{\findim}{{\rm fin.dim}}
\newcommand{\End}{{\rm End}}
\newcommand{\Ext}{{\rm Ext}}
\newcommand{\Tr}{{\rm Tr}\,}
\newcommand{\Coker}{{\rm Coker}}
\newcommand{\Ker}{{\rm Ker}}
\newcommand{\cpx}[1]{#1^{\bullet}}
\newcommand{\D}[1]{\mathscr{D}(#1)}

\newcommand{\Db}[1]{ \mathscr{D}^{\rm b}(#1)}
\newcommand{\C}[1]{\mathscr{C}(#1)}

\newcommand{\K}[1]{\mathscr{K}(#1)}

\newcommand{\Kf}[1]{\mathscr{K}^-(#1)}
\newcommand{\Kb}[1]{ \mathscr{K}^{\rm b}(#1)}
\newcommand{\modcat}[1]{#1\mbox{{\rm -mod}}}
\newcommand{\stmodcat}[1]{#1\mbox{{\rm -{\underline{mod}}}}}
\newcommand{\pmodcat}[1]{#1\mbox{{\rm -proj}}}

\newcommand{\opp}{^{\rm op}}

\newcommand{\lra}{\longrightarrow}
\newcommand{\ra}{\rightarrow}
\newcommand{\lraf}[1]{\stackrel{#1}{\lra}}
\newcommand{\raf}[1]{\stackrel{#1}{\ra}}

\begin{document}

{\Large \bf
\begin{center}
Derived equivalences from cohomological approximations, \\
and mutations of $\Phi$-Yoneda algebras

\end{center}}
\medskip

\centerline{{\bf Wei Hu}, {\bf Steffen Koenig} and {\bf Changchang
Xi$^*$}}
\bigskip

\bigskip
\begin{abstract}
In this article, a new construction of derived equivalences is
given. It relates different endomorphism rings and more generally
cohomological endomorphism rings - including higher extensions - of
objects in triangulated categories. These objects need to be
connected by certain universal maps that are cohomological
approximations and that exist in very general circumstances. The
construction turns out to be applicable in a wide variety of
situations, covering finite dimensional algebras as well as certain
infinite dimensional algebras, Frobenius categories and
$n$-Calabi-Yau categories.
\end{abstract}

\renewcommand{\thefootnote}{\alph{footnote}}
\setcounter{footnote}{-1} \footnote{ $^*$ Corresponding author.
Email: xicc@bnu.edu.cn; Fax: 0086 10 58802136; Tel.: 0086 10
58808877.}
\renewcommand{\thefootnote}{\alph{footnote}}
\setcounter{footnote}{-1} 
\renewcommand{\thefootnote}{\alph{footnote}}
\setcounter{footnote}{-1} 

\section{Introduction}

Derived equivalences have become increasingly important in
representation theory, Lie theory and geometry. Examples are ranging
from mirror symmetry over non-commutative geometry to the
Kazhdan-Lusztig conjecture and to Brou\'e's conjecture for blocks of
finite groups. In all of these situtations, and in many others,
derived equivalences are used that involve finite or infinite
dimensional algebras. Derived equivalences between algebras, or
rings, exist if and only if there exist suitable tilting complexes,
as explained quite satisfactorily by Rickard's
Morita theory for derived categories of rings
(see \cite{RickMoritaTh})).
 Derived equivalences have been
shown to preserve many significant algebraic and geometric
invariants and often to provide unexpected and useful new
connections.
\smallskip

A crucial question in this context has, however, not yet received
enough answers:

{\em How to construct derived equivalences between rings in a
general setup?}

A good answer - certainly not unique - to this question should be
general, flexible and systematic and apply to a multitude of
algebraic and geometric situations.
\smallskip

One well-developed approach is based on the theory of tilting
modules, building upon results by Happel \cite{HappelTriangle}.
Other answers use ring theoretic constructions, such as trivial
extensions \cite{Rickard3}.

The aim of this article is to provide a rather different approach.
The input of the technology developed here is a triple of objects
$(X,M,Y)$ in a triangulated category. These objects are required to
be related by certain universal maps (cohomological approximations -
a new concept introduced here, continuing approximation theory of
Auslander, Reiten and Smal\o \, \cite{arsbook}) and some
cohomological orthogonality conditions in degrees different from
zero only. The output is a derived equivalence between cohomological
endomorphism rings of $X \oplus M$ and of $M \oplus Y$.

The flexibility of the construction lies in the following features:
We enhance endomorphism rings by higher extensions to produce
cohomological endomorphism rings, broadening the classical concept
of Yoneda extension algebras. Here, we can choose a set of
cohomological degrees to define the cohomological endomorphism ring.
Choosing degree zero only gives endomorphism rings in the usual
sense - and then no orthogonality assumption is needed. Choosing all
integers, or a suitable subset thereof (satisfying an associativity
constraint), amplifies the concept of Yoneda extension algebras
$\oplus_j \Ext^j(S,S)$. There is also some flexibility in the choice
of $M$.

A special case of such a triple is given by any Auslander-Reiten
triangle $X \rightarrow M \rightarrow Y$ in a derived module
category; this already indicates generality of the construction. Our
assumptions are actually much more general and not limited to
objects in derived categories of algebras.
\smallskip

A particular feature of the derived equivalences constructed by this
method is that they also provide a very general {\em mutation
procedure}, turning one ring into another one in a systematic way.
Tilting theory has arisen as a far reaching extension of reflection
functors for quivers. Under some assumptions, but not in general, it
provides mutation procedures between two given quivers or algebras,
both of which are endomorphism rings of tilting modules; in the case
of quivers one may reflect at sink or source vertices. Mutations
similar in style also have come up in various geometric situations.
The theory of cluster categories, or more generally of Calabi-Yau
categories, has extended reflections to a mutation procedure, which
works for representations of quivers at all vertices. Such mutations
fit into the present framework as well. There is, though, a new
feature introduced by our approach: Reflection does not work in
general in derived categories (of quivers or algebras). Therefore
cluster theory passes to the cluster category, a 'quotient' of a
derived category modulo the action of some functor; endomorphism
rings are taken there. In contrast to this, the current approach
always produces equivalences on the level of derived categories, not
just of quotient categories; throughout we are considering derived
equivalences between (cohomological) endomorphism rings or quotients
thereof. In the case of quivers, this possibility of passing to
quotient algebras allows mutation at an arbitrary vertex.
\smallskip

More generality and flexibility is added by extending the concept of
'higher extensions', that is of shifted morphisms; it is possible to
replace the shift functor by any other auto-equivalence of the
ambient triangulated category. There is even a version using two
such functors.

\bigskip

The main result of this article provides a construction of derived
equivalences in a setup that is very general in several respects. In
the following explanation we start with a special case and then add
generality step by step, finally arriving at the main result.

The setup always is a triangulated category $\cal T$, which is an
$R$-category for some commutative artinian ring $R$, with identity;
so, morphism sets in $\cal T$ are $R$-modules.

\begin{enumerate}
\item
To start with, we choose any object $M$ in $\cal T$ and a triangle
$X \stackrel{\alpha}{\rightarrow} M_1 \stackrel{\beta}{\rightarrow}
Y \ra X[1]$, where $\alpha$ and $\beta$ are
$\add(M)$-approximations, that is universal maps from $X$ to objects
in $\add(M)$ or from $\add(M)$ to $Y$, respectively; in particular,
$M_1$ is in $\add(M)$. For instance, Auslander-Reiten triangles
(over algebras) provide such situations. If the triangle is induced
by an exact sequence in an abelian category, then the theorem
implies a derived equivalence between the two endomorphism rings
$\End_{\cal T}(X \oplus M_1)$ and $\End_{\cal T}(M_1 \oplus Y)$.
This can be seen as a mutation procedure relating the two
endomorphism rings. The derived equivalence has already been
established in \cite{hx2}.
\item
In the second step, recasting an idea of \cite{HX4}, endomorphism
rings are replaced by {\em cohomological endomorphism rings} in the
following sense: Higher extensions between modules $S$ and $T$ are
shifted morphisms in the derived category, $\Ext^j(S,T) \simeq
\Hom(S,T[j])$. Using Yoneda multiplication of extensions, this
defines an algebra structure on the cohomological endomorphism ring,
or generalised Yoneda algebra, $\oplus_{j \in \mathbb Z}
\Hom(S,S[j])$. When $S$ is a complex, or any object in a
triangulated category $\cal T$, negative degrees $j$ may occur. The
main theorem provides derived equivalences between such generalised
Yoneda algebras. The construction works, however, not only for these
Yoneda algebras, but also for 'perforated' ones in the following
sense: Choose a subset $\Phi \subset \mathbb Z$. Then, under some
associativity constraint requiring $\Phi$ to be 'admissible' (see
Subsection \ref{2.3}), the space $\oplus_{j \in \Phi} \Hom(S,S[j])$
is an associative algebra, that in general is neither a subalgebra
nor a quotient algebra of the Yoneda algebra $\oplus_{j \in \mathbb
Z} \Hom(S,S[j])$. This algebra is called a {\em $\Phi$-Yoneda
algebra} or a {\em $\Phi$-perforated Yoneda algebra}. We will use
the notation $\Ex{\cal T}{\Phi}{Z}$ for the algebra $\oplus_{j \in
\Phi} \Hom(Z,Z[j])$, where $Z$ is any object in $\cal T$.
\\
The assumptions of the first step get modified by using
cohomological approximations, in the degrees specified by $\Phi$,
instead of approximations in degree zero only. Auslander-Reiten
triangles still satisfy these properties. Adding higher extensions
requires also to add an orthogonality assumption without which the
result would be wrong: Assume $\Hom(M,X[j]) = 0 = \Hom(Y,M[j])$ for
all $j \in \Phi, j \neq 0$. For the sake of exposition also assume
for a moment that the above triangle $X
\stackrel{\alpha}{\rightarrow} M_1 \stackrel{\beta}{\rightarrow} Y
\ra X[1]$ is in a derived module category and it is induced from an
exact sequence with corresponding properties. Then there are derived
equivalences between $\Phi$-Yoneda algebras $\Db{\Ex{\cal
T}{\Phi}{X\oplus M}} \simeq \Db{\Ex{\cal T}{\Phi}{M\oplus Y}}$.
\item
This result needs to be modified, if the triangle is not induced by
an exact sequence any more. Then some annihilators have to be
factored out of the degree zero parts of the cohomological
endomorphism rings, and the derived equivalences are connecting the
quotient algebras $\Ex{\cal T}{\Phi}{X\oplus M}/I$ and $\Ex{\cal
T}{\Phi}{M\oplus Y}/J$. Here, the ideals $I$ and $J$ can be
described as follows: Let $\Gamma_0=\End_{\cal T}(M \oplus Y)$ and
$e$ the idempotent element in $\Gamma_0$ corresponding to the direct
summand $M$. Then $J$ is the submodule of the left $\Gamma_0$-module
$\Gamma_0e\Gamma_0$, which is maximal with respect to $eJ=0$. Let
$\Lambda_0=\End_{\cal T}(X\oplus M)$, and $f$ the idempotent in
$\Lambda_0$ corresponding to the direct summand $M$. Then $I$ is the
submodule of the right $\Lambda_0$-module $\Lambda_0f\Lambda_0$,
which is maximal with respect to $If=0$.

Another, equivalent, description of $I$ and $J$ is that $I$ consists
of all elements $(x_i)_{i\in \Phi}\in \Ex{\cal T}{\Phi}{X\oplus M}$
such that $x_i=0$ for $0\ne i\in \Phi$ and $x_0$ factorises through
$\add(M)$ and $x_0\widetilde{\alpha}=0$, and $J$ consists of all
elements $(y_i)_{i\in \Phi}\in \Ex{\cal T}{\Phi}{M\oplus Y}$ such
that $y_i=0$ for $0\ne i\in \Phi$ and $y_0$ factorises through
$\add(M)$ and $\bar{\beta}y_0=0$, where $\widetilde{\alpha}$ is the
diagonal morphism diag$(\alpha,1): X\oplus M\ra M_1\oplus M$, and
$\bar{\beta}$ is the skew-diagonal morphism skewdiag$(1, \beta):
M_1\oplus M\ra M\oplus Y$.
\item
The fourth level of generalisation allows to replace the shift
functor by any auto-equivalence of the triangulated category
$\mathcal T$, thus providing a new and versatile meaning of 'higher
extensions' in terms of morphisms with one variable shifted by
powers of the auto-equivalence. The additional datum $F$ gets
mentioned, when necessary, in the notation as an additional
superscript, as in $\Ex{\cal T}{F,\Phi}{Z}$.
\end{enumerate}

In this general form, the main theorem is as follows:

\begin{Theo}
Let $\Phi$ be an admissible subset of $\mathbb Z$, and let $\cal T$
be a triangulated $R$-category and $M$ an object in $\cal T$. Assume
that $F$ is an invertible triangle functor from $\cal T$ to itself.
Suppose that $$X\stackrel{\alpha}{\lra} M_1\stackrel{\beta}{\lra} Y
\stackrel{w}{\lra} X[1]$$ is a triangle in $\mathcal T$ such that

\smallskip
$(1)$ The morphism $\alpha$ is a left $(\add(M), F,
\Phi)$-approximation of $X$ and $\beta$ is a right $(\add(M), F,
{-\Phi})$-approximation of $Y$,

$(2)$ $\Hom_{\cal T}(M,F^iX)=0= \Hom_{\cal T}(F^{-i}Y, M)$ for all
$0\ne i\in \Phi$.

\smallskip
\noindent Then $\Ex{\cal T}{F,\Phi}{X\oplus M}/I$ and $\Ex{\cal
T}{F,\Phi}{M\oplus Y}/J$ are derived equivalent, where $I$ and $J$
are the above ideals of the $\Phi$-Yoneda algebras $\Ex{\cal
T}{F,\Phi}{X\oplus M}$ and $\Ex{\cal T}{F,\Phi}{M\oplus Y}$,
contained in $\End_{\cal T}(X\oplus M)$ and $\End_{\cal T}(M\oplus
Y)$, respectively. \label{thm1}
\end{Theo}

A fifth level of generalisation, using two functors $F$ and $G$,
will be discussed in the Appendix. A further generalisation of some
results in this paper to $n$-angulated categories introduced in
\cite{GKO} will be considered in \cite{ChYP}.
\medskip

The second level of generality, where $F$ is the shift functor and
both $I$ and $J$ are zero, is already widely applicable. This case
happens frequently for the derived category $\Db{A}$ of an
$R$-algebra $A$.

\begin{Koro}
Let $\Phi$ be an admissible subset of $\mathbb N$, and let $A$ be an
$R$-algebra and $M$ an $A$-module. If $0\ra X\stackrel{\alpha}{\lra}
M_1\stackrel{\beta}{\lra} Y \ra 0$ is an exact sequence in
$\modcat{A}$ such that $\alpha$ is a left $(\add(M),
\Phi)$-approximation of $X$ and $\beta$ is a right $(\add(M),
-\Phi)$-approximation of $Y$ in $\Db{A}$,  and that
$\Ext^i_{A}(M,X)= 0 = \Ext^i_{A}(Y, M)$ for all $0\ne i\in \Phi$,
then the $\Phi$-Yoneda algebras $\Ex{A}{\Phi}{X\oplus M}$ and
$\Ex{A}{\Phi}{M\oplus Y}$ are derived equivalent. \label{cor2}
\end{Koro}

These results partly generalise some results of \cite{hx2}.
\bigskip

The setup here, and the main result, covers, combines and extends
several classical concepts:

Auslander algebras - endomorphism rings of direct sums of 'all'
modules of an algebra of finite representation type - are the
ingredients of the celebrated Auslander correspondence,
characterising finite representation type via homological
dimensions. Auslander algebras of derived equivalent algebras are,
in general, not derived equivalent; positive results in this
direction - for self-injective algebras of finite representation
type - previously have been obtained in \cite{HX4}. In the current
approach new results can be obtained by appropriate choices of $X
\oplus M$.

Another intensively studied class of algebras is that of Yoneda
algebras, that is, algebras of self-extensions of a semisimple
module, or more generally of any module. Apparently, the
constructions in Corollary \ref{cor2} and in \cite{HX4} provide the
first general class of derived equivalences for Yoneda algebras.
Perforated Yoneda algebras first have been defined in \cite{HX4},
under the name $\Phi$-Auslander-Yoneda algebras. The approach
developed there has been based on the existence of particular kinds
of derived equivalences for algebras, which then have been used to
construct derived equivalences for perforated Yoneda algebras.

The main novelty of the present approach is the systematic use of
cohomological data, such as cohomological approximations and
perforated Yoneda algebras. This relates smoothly with a wide
variety of concepts, such as Auslander-Reiten sequences and
triangles, dominant dimension, Calabi-Yau categories and Frobenius
categories.
\bigskip

The article is organised as follows. In Section \ref{sect2}, we
first fix notation, and then recall definitions and basic results on
derived equivalences as well as on admissible sets and perforated
Yoneda algebras. Also, we extend the notion of $\cal
D$-approximation to what we call cohomological ${\cal
D}$-approximation with respect to $(F,\Phi)$, where $F$ is a functor
and $\Phi$ is a subset of $\mathbb N$. In Section \ref{sect3}, the
main result, Theorem \ref{thm1}, is proven and various easier to
access situations are described, for which the assumptions of
Theorem \ref{thm1} are satisfied. Section \ref{sect4} explains how
Theorem \ref{thm1} applies to a variety of situations: derived
categories of Artin algebras, Frobenius categories and Calabi-Yau
categories. Also, the connection to the concept of dominant
dimension is explained. In Section \ref{sect5}, two examples are
given to illustrate the results and to show the necessity of some
assumptions in Theorem \ref{thm1}. In the Appendix, a more general
formulation of Theorem \ref{thm1} is stated, which involves two
functors, in order to add more flexibility with a view to potential
future applications.

\bigskip
The authors are grateful to Rundong Zheng and Yiping Chen at BNU for carefully reading the first version of the manuscript. The corresponding
author C.C. Xi thanks NSFC for partial support. W. Hu is
grateful to the Alexander von Humboldt Foundation for a Humboldt
fellowship. Much of this work has been done during visits of
Xi and Hu to the Mathematisches Institut, Universit\"at zu K\"oln,
in 2010.

\section{Preliminaries \label{sect2}}

In this section, we shall recall basic definitions and facts which
will be needed in the proofs later on.

\subsection{Conventions}

Throughout this paper, $R$ is a fixed commutative artinian ring with
identity. Given an $R$-algebra $A$, by an $A$-module we mean a
unitary left $A$-module; the category of all (respectively, finitely
generated) $A$-modules is denoted by $A$-Mod (respectively,
$\modcat{A}$), the full subcategory of $A$-Mod consisting of all
(respectively, finitely generated) projective modules is denoted by
$A$-Proj (respectively, $\pmodcat{A}$). There is a similar notation
for right $A$-modules. The stable module category $\stmodcat{A}$ of
$A$ is, by definition, the quotient category of $\modcat{A}$ modulo
the ideal generated by homomorphisms factorising through projective
modules in $\pmodcat{A}$. An equivalence between the stable module
categories of two algebras is called a {\em stable equivalence}.

An $R$-algebra $A$ is called an \emph{Artin $R$-algebra} if $A$ is
finitely generated as an $R$-module. For an Artin $R$-algebra $A$,
we denote by $D$ the usual duality on $\modcat{A}$, and by $\nu_A$
the Nakayama functor $D\Hom_A(-, {}_AA): A\mbox{-proj}\ra
A\mbox{-inj}$. For an $A$-module $M$, we denote the first syzygy of
$M$ by $\Omega_A(M)$, and call $\Omega_A$ the \emph{Heller loop
operator} of $A$. The transpose of $M$, which is an $A\opp$-module,
is denoted by $\Tr(M)$.

Let $\cal C$ be an additive $R$-category, that is, $\cal C$ is an
additive category in which the set of morphisms between two objects
in $\cal C$ is an $R$-module, and the composition of morphisms in
$\cal C$ is $R$-bilinear. For an object $X$ in $\mathcal{C}$, we
denote by $\add(X)$ the full subcategory of $\cal C$ consisting of
all direct summands of finite direct sums of copies of $X$. An
object $X$ in $\cal C$ is called an \emph{additive generator} for
$\cal C$ if ${\cal C}$ = add$(X)$. For two morphisms $f:X\rightarrow
Y$ and $g:Y\rightarrow Z$ in $\cal C$, we write $fg$ for their
composition. For two functors $F:\mathcal{C}\rightarrow \mathcal{D}$
and $G:\mathcal{D}\rightarrow\mathcal{E}$ however, we write $GF$ for
the composition instead of $FG$.

If $f: X\ra Y$ is a map between two sets $X$ and $Y$, we denote the
image of $f$ by Im$(f)$. Moreover, if $f$ is a homomorphism between
two abelian groups, we denote the kernel and cokernel of $f$ by
$\Ker(f)$ and $\Coker(f)$, respectively.

Recall that a functor $F:{\cal C}\ra {\cal D}$ is invertible if
there is a functor $G: {\cal D}\ra {\cal C}$ such that $GF=\id_{\cal
C}$ and $FG=\id_{\cal D}$. In this case we write $F^{-1}$ for $G$.
If $\cal C$ = $\cal D$ and if $F$ is invertible, we say that $F$ is
an auto-equivalence. If $F$ is a functor from $\cal C$ to $\cal C$,
then we write $F^0=\id_{\cal C}$, and $F^{-i} = (F^{-1})^i$ for $i>
0$ if $F^{-1}$ exists, and $F^{-i}=0$ otherwise.

Let $\cal T$ be a triangulated $R$-category with a shift functor
[1]. For two objects $X$ and $Y$ in $\cal T$, we write sometimes
$\Ext^i_{\cal T}(X,Y)$ for $\Hom_{\cal T}(X,Y[i])$. Let $\Phi$ be a
subset of $\mathbb Z$. An object $M$ (or a full subcategory $\cal
M$) of ${\cal T}$ is called \emph{$\Phi$-self-orthogonal} provided
that $\Ext^i_{\cal T}(M, M)=0$ (or $\Ext^i_{\cal T}({\cal M},{\cal
M})=0$ ) for all $0\neq i\in\Phi$, where $\Ext^i_{\cal T}({\cal
M},{\cal M})=0$ means that $\Ext^i_{\cal T}(X,Y)=0$ for all $X,Y\in
{\cal M}$. In case $\Phi =\mathbb Z$, we say that $M$ is
\emph{self-orthogonal}. For $\Phi = \{0, 1, \cdots, n\}$, we say
that $M$ is $n$-self-orthogonal, which is sometimes, perhaps less
suggestively, referred to as $n$-rigid.
%It seems that ``self-orthogonal" is more suggestive
%than the word ``rigid" for an arbitrary $\Phi$.

Replacing the shift functor by a triangle auto-equivalence $F$, one
may also define the notion of $(F,\Phi)$-self-orthogonality, but we
refrain from introducing this notion here.

\subsection{Derived equivalences}

Let $\cal C$ be an additive $R$-category.

By a complex $\cpx{X}$ over $\cal C$ we mean a sequence of morphisms
$d_X^{i}$ between objects $X^i$ in $\cal C$: $ \cdots \rightarrow
X^i\stackrel{d_X^i}{\lra}
X^{i+1}\stackrel{d_X^{i+1}}{\lra}X^{i+2}\rightarrow\cdots, $ such
that $d_X^id_X^{i+1}=0$ for all $i \in {\mathbb Z}$; we write
$\cpx{X}=(X^i, d_X^i)$. For a complex $\cpx{X}$, the {\em brutal
truncation} $\sigma_{<i}\cpx{X}$ of $\cpx{X}$ is a quotient complex
of $\cpx{X}$ such that $(\sigma_{<i}\cpx{X})^k$ is $X^k$ for all
$k<i$ and zero otherwise. Similarly, we define $\sigma_{\geqslant
i}\cpx{X}$. For a fixed $n\in {\mathbb Z}$, we denote by
$\cpx{X}[n]$ the complex obtained from $\cpx{X}$ by shifting degree
by $n$, that is, $(\cpx{X}[n])^0=X^n$.

The category of all complexes over $\cal C$ with chain maps is
denoted by $\C{\mathcal C}$. The homotopy category of complexes over
$\mathcal{C}$ is denoted by $\K{\mathcal C}$. When $\cal C$ is an
abelian category, the derived category of complexes over $\cal C$ is
denoted by $\D{\cal C}$. The full subcategories of $\K{\cal C}$ and
$\D{\cal C}$ consisting of bounded complexes over $\mathcal{C}$ are
denoted by $\Kb{\mathcal C}$ and $\Db{\mathcal C}$, respectively. As
usual, for an algebra $A$, we simply write $\C{A}$ for
$\C{\modcat{A}}$, $\K{A}$ for $\K{\modcat{A}}$ and $\Kb{A}$ for
$\Kb{\modcat{A}}$. Similarly, we write $\D{A}$ and $\Db{A}$ for
$\D{\modcat{A}}$ and $\Db{\modcat{A}}$, respectively.

For an $R$-algebra $A$, the categories $\K{A}$ and $\D{A}$ are
triangulated $R$-categories. For basic results on triangulated
categories, we refer the reader to \cite{HappelTriangle} and
\cite{neeman}.

The following result, due to Rickard (see \cite[Theorem
6.4]{RickMoritaTh}) by a direct approach, and to Keller by working
in the more general setup of differential graded algebras, 
%(see, for instance, \cite{kellerdc}), 
is fundamental in the investigation of derived equivalences.

\begin{Theo}{\rm \cite{RickMoritaTh}}\label{rickard}
Let $\Lambda$ and $\Gamma$ be two rings. The following conditions
are equivalent:

$(a)$ $\Kf{\Lambda\emph{-Proj}}$ and $\Kf{\Gamma\emph{-Proj}}$ are
equivalent as triangulated categories;

$(b)$ $\Db{\Lambda\emph{-Mod}}$ and $\Db{\Gamma\emph{-Mod}}$ are
equivalent as triangulated categories;

$(c)$ $\Kb{\Lambda\emph{-Proj}}$ and $\Kb{\Gamma\emph{-Proj}}$ are
equivalent as triangulated categories;

$(d)$ $\Kb{\Lambda\emph{-proj}}$ and $\Kb{\Gamma\emph{-proj}}$ are
equivalent as triangulated categories;

$(e)$ $\Gamma$ is isomorphic to
$\End_{\Kb{\Lambda\emph{-proj}}}(\cpx{T})$, where $\cpx{T}$ is a
complex in $\Kb{\Lambda\emph{-proj}}$ satisfying:

\qquad $(1)$ $\cpx{T}$ is self-orthogonal, that is,
$\Hom_{\Kb{\pmodcat{\Lambda}}}(\cpx{T},\cpx{T}[i])=0$ for all $i\ne
0$,

\qquad $(2)$ $\add(\cpx{T})$ generates $\Kb{\pmodcat{\Lambda}}$ as a
triangulated category. \label{derequ}
\end{Theo}

Two rings $\Lambda$ and $\Gamma$ are called \emph{derived
equivalent} if the above conditions (a)-(e) are satisfied. A complex
$\cpx{T}$ in $\Kb{\pmodcat{\Lambda}}$ as above is called a
\emph{tilting complex} over $\Lambda$.

For Artin algebras, the above equivalent conditions can be
reformulated in terms of finitely generated modules: Two Artin
$R$-algebras $A$ and $B$ are said to be \emph{derived equivalent} if
their derived categories $\Db{A}$ and $\Db{B}$ are equivalent as
triangulated categories. In this case, there is a tilting complex
$\cpx{T}$ in $\Kb{\pmodcat{A}}$ such that $B\simeq
\End_{\Kb{A}}(\cpx{T})$.

\subsection{Admissible subsets and $\Phi$-Yoneda algebras\label{2.3}}

Let ${\mathbb N}=\{0,1, 2, \cdots\}$ be the set of natural numbers,
and let $\mathbb Z$ be the set of all integers. For a natural number
$n$ or infinity, let ${\mathbb N}_n :=\{i\in {\mathbb N}\mid 0\le i<
n+1\}$.

Recall from \cite{HX4} that a subset $\Phi$ of $\mathbb{Z}$
containing $0$ is called an {\em admissible subset} of $\mathbb{Z}$
if the following condition is satisfied:

{\it If $i, j$ and $k$ are in $\Phi$ such that $i+j+k\in\Phi$, then
$i+j\in\Phi$ if and only if $j+k\in\Phi$.}

\medskip
Any subset $\{0,i,j\}$ of $\mathbb N$ is an admissible subset of
$\mathbb Z$. Moreover, for any subset $\Phi$ of $\mathbb N$
containing zero and for any positive integer $m\ge 3$, the set
$\{x^m\mid x\in \Phi\}$ is admissible in $\mathbb Z$ (for more
examples, see \cite{HX4}). Nevertheless, not every subset of
$\mathbb N$ containing zero is admissible, for instance,
$\{0,1,2,4\}$ is not admissible. In fact, this is the 'smallest'
non-admissible subset of $\mathbb N$.

Admissible sets were used to define $\Phi$-Yoneda algebras in
\cite{HX4}, under the name of '$\Phi$-Auslander-Yoneda algebras'.
The formulation there works more generally for monoid graded
algebras. For our purpose in this paper, we restrict to the case of
an object in a triangulated category.

Let $\Phi$ be an admissible subset of $\mathbb{Z}$, and let
$\mathcal T$ be a triangulated $R$-category with a shift functor
[1]. Suppose that $F$ is a triangle functor from $\cal T$ to $\cal
T$. Recall that we put $F^i=0$ for $i<0$ if $F^{-1}$ does not exist.

Let $E_{\mathcal T}^{F,\Phi}(-,-)$ be the bi-functor
$$\bigoplus_{i\in\Phi}\, \Hom_{\cal T}(-, F^i-): {\mathcal T}\times{\mathcal T}\lra R\mbox{\rm -Mod},$$
$$ (X,Y)\mapsto E_{\mathcal T}^{F,\Phi}(X,Y):= \bigoplus_{i\in\Phi}\, \Hom_{\mathcal T}(X, F^iY),$$
$$ X\lraf{f} X' \mapsto \bigoplus_{i\in \Phi} \Hom_{\cal T}(f,F^iY),\qquad Y\lraf{g} Y'\mapsto \bigoplus_{i\in\Phi}\Hom_{\cal T}(X,F^ig).$$
Suppose that $X, Y$ and $Z$ are objects in $\mathcal T$. Let
$(f_i)_{i\in \Phi}\in\mbox{E}_{\cal T}^{F,\Phi}(X,Y)$ and
$(g_i)_{i\in \Phi}\in\mbox{E}_{\cal T}^{F,\Phi}(Y,Z)$. We define a
composition as follows:

$$\mbox{E}_{\cal T}^{F,\Phi}(X,Y)\times\mbox{E}_{\cal T}^{F,\Phi}(Y,Z)\lra\mbox{E}_{\cal T}^{F,\Phi}(X,Z)$$
$$ \big((f_i)_{i\in \Phi},(g_i)_{i\in \Phi}\big)\mapsto \big(\sum_{{{u, v\in\Phi}\atop {u+v=i}}}f_u(F^ug_v)\big)_{i\in \Phi}.$$

Since $\Phi$ is admissible, this composition is associative. Thus
$\mbox{E}_{\cal T}^{F,\Phi}(X,X)$ is an $R$-algebra. It is called
the \emph{ $\Phi$-Yoneda algebra} or, when $\Phi$ is fixed, the
\emph{perforated Yoneda algebra} of $X$ with respect to $F$. Then
$\mbox{E}_{\cal T}^{F,\Phi}(X,Y)$ is a left $\mbox{E}_{\cal
T}^{F,\Phi}(X,X)$-module. When $\Phi={\mathbb N}$, the algebra
$\E{\cal T}{F,\Phi}{X}{X}$ is the orbit algebra of $X$ under $F$
(see \cite{DL}).

For convenience we write $\Ex{\cal T}{F,\Phi}{X}$ for $\E{\cal
T}{F,\Phi}{X}{X}$.  In case ${\cal T} = \Db{A}$ where $A$ is a ring
with identity, we write $\E{A}{F,\Phi}{X}{Y}$ for
$\E{\Db{A}}{F,\Phi}{X}{Y}$, and $\Ex{A}{F,\Phi}{X}$ for
$\Ex{\Db{A}}{F,\Phi}{X}$.

When $F$ coincides with the shift functor, we omit the upper index
$F$, and call $\Ex{\cal T}{\Phi}{X}$ the $\Phi$-Yoneda algebra of
$X$, without referring to the shift functor. This is the algebra
introduced in \cite{HX4} and there called an Auslander-Yoneda
algebra.

The following lemma is essentially taken from \cite[Lemma 3.5]{HX4},
where a variation of it appears. The proof given there carries over
to the present situation.

\begin{Lem}
Let $\cal T$ be a triangulated $R$-category with a triangle
endo-functor $F$, and let $U$ be an object in $\cal T$. Suppose that
$U_1$, $U_2$ and $U_3$ are in $\add(U)$, and that $\Phi$ is an
admissible subset of $\mathbb Z$. Then,

\medskip
$(1)$ there is a natural isomorphism
$$\mu: \E{\cal T}{F,\Phi}{U_1}{U_2}\lra
\Hom_{\Ex{\cal T}{F,\Phi}{U}}(\E{\cal T}{F,\Phi}{U}{U_1}, \E{\cal
T}{F,\Phi}{U}{U_2}),$$ which sends $x\in\E{\cal
T}{F,\Phi}{U_1}{U_2}$ to the morphism $a \mapsto ax$ for $a\in
\E{\cal T}{F,\Phi}{U}{U_1}$. Moreover, if $x\in\E{\cal
T}{F,\Phi}{U_1}{U_2}$ and $y\in\E{\cal T}{F,\Phi}{U_2}{U_3}$, then
$\mu(xy)=\mu(x)\mu(y)$.

\smallskip
$(2)$ The functor $\E{\cal T}{F,\Phi}{U}{-}: \add(U)\lra
\pmodcat{\Ex{\cal T}{F,\Phi}{U}}$ is faithful.

\smallskip
$(3)$ If $\Hom_{\cal T}(U_1,F^iU_2)=0$ for all $i\in\Phi\setminus
\{0\}$, then the functor $\E{\cal T}{F,\Phi}{U}{-}$ induces an
isomorphism of $R$-modules:
$$\E{\cal T}{F,\Phi}{U}{-}: \Hom_{\cal T}(U_1, U_2)\lra
\Hom_{\Ex{\cal T}{F,\Phi}{U}}(\E{\cal T}{F,\Phi}{U}{U_1}, \E{\cal
T}{F,\Phi}{U}{U_2}).$$ \label{Extprop}
\end{Lem}
The properties described in Lemma \ref{Extprop} will be frequently
used in the proofs below.

The class of $\Phi$-Yoneda algebras with respect to a functor
includes a large class of algebras, for example the following:

({\bf a}) The endomorphism algebra of a module, in particular, the
Auslander algebras of representation-finite algebras. Here we choose
$\Phi=\{0\}$.

({\bf b}) The generalised Yoneda algebra of a module if we take
$\Phi=\mathbb N$. This includes the preprojective algebras (see
\cite{DL}) and the Hochschild cohomology rings of given algebras.
Choosing $\Phi=2{\mathbb N}$, we get for instance the even
Hochschild cohomology rings of algebras.

({\bf c}) Certain trivial extensions: For an Artin algebra $A$ and
an $A$-module $M$ we choose $\Phi=\{0,i\}$ for $i\ge 1$ an arbitrary
natural number. Then $\Ex{A}{\Phi}{M}$ is the trivial extension of
$\End_A(M)$ by the bimodule $\Ext^i_A(M,M)$. Such rings appear
naturallly in the (bounded) derived category $\Db{\mathbb X}$ of
coherent sheaves of a smooth projective variety $\mathbb X$ over
$\mathbb C$. Indeed, if $X$ is a $d$-spherical object in
$\Db{\mathbb X}$, then its cohomological ring
$\cpx{\End}_{\Db{\mathbb X}}(X)$ is $\Ex{\Db{\mathbb
X}}{\{0,d\}}{X}$, this is a graded ring isomorphic to ${\mathbb
C}[t]/(t^2)$ with $t$ of degree $d$. For further information on
spherical objects, we refer the reader to \cite[Section 3c]{st}.

In general, if $\Phi=\{0, a_1, \cdots, a_n\}\subseteq {\mathbb N}$
such that $a_i>2a_{i-1}$ for $i=2, \cdots, n$, then
$\Ex{A}{\Phi}{X}$ is the trivial extension of $\End_A(X)$ by the
bimodule $\displaystyle\bigoplus_{0\ne i\in \Phi}\Ext^i_A(X,X)$.
Note that $\Phi =\{0\}\cup \{2n+1\mid n\in{\mathbb N}\}$ is
admissible. In this case, we also get a trivial extension.

({\bf d}) The polynomial ring $R[t]$: If we take $\Phi=m{\mathbb N}$
for $m\ge 1$, then the perforated Yoneda algebra
$\Ex{R[x]/(x^2)}{\Phi}{R}$ is isomorphic to $R[t^m]$ with $t$ a
variable. If $\Phi=\{0, 1, \cdots, n\}$, then
$\Ex{R[x]/(x^2)}{\Phi}{R}\simeq R[t]/(t^n)$.

\subsection{ $\cal D$-split sequences and cohomological
$\cal D$-approximations\label{subsect2.4}}

 $\cal D$-split sequences have been defined in \cite{hx2} in
the context of constructing derived equivalences between certain
endomorphism algebras. Let us recall the definition and a result in
\cite{hx2}.

Let $\cal C$ be an additive category and $\cal D$ a full subcategory
of $\cal C$. A sequence
$$X\stackrel{f}{\longrightarrow}M\stackrel{g}{\longrightarrow}Y$$ in $\cal C$
is called an \emph{$\cal D$-split sequence} if

\smallskip
$(1)$ $M\in {\cal D}$,

$(2)$ $f$ is a left $\cal D$-approximation of $X$, and $g$ is a
right $\cal D$-approximation of $Y$, and

$(3)$ $f$ is a kernel of $g$, and $g$ is a cokernel of $f$.

\medskip
Typical examples of $\cal D$-split sequences are Auslander-Reiten
sequences. Every $\mathcal D$-split sequence provides a derived
equivalence (see \cite[Theorem 1.1]{hx2}). Here are some details,
for later reference.

\begin{Theo} {\rm \cite{hx2}}
Let $\cal C$ be an additive category, and $M$ an object in $\cal C$.
Suppose that
$$X\longrightarrow M'\longrightarrow Y
$$ is an  \emph{add}$(M)$-split sequence in
$\cal C$. Then the endomorphism ring $\End_{\cal C}(M\oplus X)$ of
$M\oplus X$ is derived-equivalent to the endomorphism ring
$\End_{\cal C}(M\oplus Y)$ of $M\oplus Y $ via a tilting module of
projective dimension at most $1$. \label{deadss}
\end{Theo}

Now, the question arises whether Theorem \ref{deadss} can be
extended to $\Phi$-Yoneda algebras. The second example in the final
section demonstrates that this is no longer true if we just replace
the endomorphism algebras in Theorem \ref{deadss} by $\Phi$-Yoneda
algebras. Nevertheless, we shall show that under certain
orthogonality conditions, there still is a positive answer. This
will be discussed in detail in the next section.

The condition (3) of a $\cal D$-split sequence are substitutes in
this general setup for requiring the short exact sequence to be
exact. Since triangles in triangulated categories are replacements
of short exact sequences, we may reformulate the notion of $\cal
D$-split sequences in the following sense for triangulated
categories.

Let $\cal T$ be a triangulated category with a shift functor [1],
and let $\cal D$ be a full additive subcategory of $\cal T$. A
triangle
$$ X\stackrel{\alpha}{\lra} M'\stackrel{\beta}{\lra} Y\lra X[1]$$
in $\cal T$ is called a \emph{$\cal D$-split triangle} if $M'\in
{\cal D}$, the map $\alpha$ is a left $\cal D$-approximation of $X$
and the map $\beta$ is a right $\cal D$-approximation of $Y$.

Thus, for an Artin $R$-algebra $A$, every $\cal D$-split sequence in
$A$-mod extends to a $\cal D$-split triangle in $\Db{A}$.
%; here, to guarantee existence, $\cal D$ has to be chosen appropriately.

\medskip
Next, we introduce the left and right cohomological $\cal
D$-approximations with respect to $(F,\Phi)$, which generalise the
notions of left and right $\cal D$-approximations, respectively.

Suppose that $\cal C$ is a category with an endo-functor $F:{\cal
C}\ra {\cal C}$. Let $\cal D$ be a full subcategory of $\cal C$, and
let $\Phi$ be a non-empty subset of $\mathbb N$. If $F$ has an
inverse, then $\Phi$ may be chosen to be a subset of $\mathbb Z$.
Suppose that $X$ is an object of $\cal C$. A morphism $f: X\ra D$ in
$\cal C$ is called a \emph{left cohomological $\cal
D$-approximation} of $X$ with respect to $(F,\Phi)$ (or shortly, a
left $({\cal D},F, \Phi)$-approximation of $X$) if $D\in {\cal D}$,
and for any morphism $g: X\ra F^i(D')$ with $D'\in {\cal D}$ and
$i\in \Phi$, there is a morphism $g': D\ra F^i(D')$ such that
$g=fg'$. Here $F^0=\id_{\cal C}$. Similarly, we have the notion of a
right $({\cal D}, F, \Phi)$-approximation of $X$ in $\cal T$, that
is, a morphism $f: D\ra X$ with $D$ in $\cal D$ is called a right
$({\cal D}, F, \Phi)$-approximation of $X$ if, for any $i\in \Phi$
and any morphism $g: F^iD'\ra X$ with $D'$ in $\cal D$, there is a
morphism $g': F^iD'\ra D$ such that $g=g'f$.

Note that if $F=\id_{\cal C}$ and $\Phi=\{0\}$, then we get the
original notion of approximations in the sense of Auslander and
Smal\o. (In ring theory, such approximations are called pre-envelope
and precover, respectively). Moreover, if $0\in \Phi$, then every
left $({\cal D}, F, \Phi)$-approximation of $X$ is also a left $\cal
D $-approximation of $X$, and every right $({\cal D}, F,
\Phi)$-approximation of $X$ is also a right $\cal D$-approximation
of $X$.

If $F=[1]$ and $\cal T$ = $\Db{A}$ for an Artin algebra $A$, then
$\Hom_{\cal T}(X,F^iY)\simeq\Ext^i_A(X,Y)$ for all $X, Y\in A$-mod
and all $i\ge 0$. For this reason, a $({\cal D},
F,\Phi)$-approximation has been called a \emph{cohomological}
approximation.

In this paper, we are mainly interested in the case where $\cal C$
is a triangulated $R$-category $\cal T$ with an endo-functor $F$,
and $\cal D$ is a full subcategory of $\cal T$. Thus, a morphism $f:
X\ra D$ with $D\in {\cal D}$ and $X\in {\cal T}$ is a left $({\cal
D}, F, \Phi)$-approximation of $X$ if and only if the canonical map
$\E{\cal T}{F,\Phi}{f}{D'}: \E{\cal T}{F,\Phi}{D}{D'}\ra \E{\cal
T}{F,\Phi}{X}{D'}$, defined by $(x_i)_{\in \Phi}\mapsto (fx_i)_{i\in
\Phi}$, is surjective for all $D'\in {\cal D}$. Similarly, a
morphism $g: D\ra X$ with $D\in {\cal D}$ and $X\in {\cal T}$ is a
right $({\cal D}, F, \Phi)$-approximation of $X$ if and only if the
canonical map $\Hom_{\cal T}(F^jD',g): \Hom_{\cal T}(F^jD', D)\ra
\Hom_{\cal T}(F^jD', X)$ is surjective for every $D'\in {\cal D}$
and $j\in \Phi$. If, moreover, $F$ is an invertible triangle
functor, then a morphism $g: D\ra X$ with $D\in {\cal D}$ and $X\in
{\cal T}$ is a right $({\cal D}, F, \Phi)$-approximation of $X$ if
and only if the canonical map $\E{\cal T}{F,-\Phi}{D'}{g}: \E{\cal
T}{F,-\Phi}{D'}{D}\ra \E{\cal T}{F, -\Phi}{D'}{X}$ is surjective for
all $D'\in {\cal D}$. Note that here we need the minus sign for
$\Phi$ and that $F^{-1}$ exists.

%To avoid the minus sign and the existence of the inverse of $F$, we
%may introduce another type of approximation. A morphism $\beta: D\ra
%X$ with $D\in \cal D$ is called a right $({\cal
%D},F\Phi)$-approximation of $X$ provided that every morphism $g:
%D'\ra F^iX$, with $D'\in {\cal D}$ and $i\in \Phi$, factorises
%through $F^i\beta: F^iD\ra F^iX$. If $F$ is invertible, then $\beta$
%is a right $({\cal D},F, -\Phi)$-approximation of $X$ if and only if
%$\beta$ is a right $({\cal D}, F\Phi)$-approximation of $X$. If $F$
%is not invertible, then a right $({\cal D},F, -\Phi)$-approximation
%may be different from a right $({\cal D}, F\Phi)$-approximation.

If $F$ coincides with the shift functor [1], we simply speak of
$({\cal D},\Phi)$-approximations, without mentioning $F$.

Note that if $\Phi$ contains zero and if $\Hom_{\cal T}(X,F^iD')=0$
for all $0\ne i\in \Phi$ and $D'\in{\cal D}$, then $f$ is a left
$({\cal D},F, \Phi)$-approximation of $X$ if and only if $f$ is a
left $\cal D$-approximation of $X$. A dual statement is also true
for a right $({\cal D}, F, \Phi)$-approximation of $X$.

%Clearly, one may use $({\cal D}, F,{\Phi})$-approximations to define
%the notion of an  $({\cal D}, F, \Phi)$-split triangle in
%$\cal T$. Here we will not state it in detail.

Here is a source of examples of $({\cal D}, \Phi)$-approximations.
Suppose that $\cal T$ = $\Db{A}$ for $A$ an Artin $R$-algebra and
that $\Phi$ is a subset of $\mathbb Z$. Let $
X\stackrel{\alpha}{\lra}M\stackrel{\beta}{\lra}Y\ra X[1]$ be an
Auslander-Reiten triangle in $\cal T$. If neither $X$ nor $Y$
belongs to $\add(M[i])$ for every $0\ne i\in \Phi$, then $\alpha$ is
a left $(\add(M),\Phi)$-approximation of $X$, and  $\beta$ is a
right $(\add(M),\Phi)$-approximation of $Y$.

Finally, we note the difference of a left $({\cal D}, F,
\Phi)$-approximation of $X$ from a left $\big(\cup_{i\in
\Phi}F^i{\cal D}\big)$-approximation of $X$ in the sense of
Auslander and Smal\o, where $\cup_{i\in \Phi}F^i{\cal D}$ is the
full subcategory of $\cal T$ with all objects in $F^i{\cal D}$ for
all $i\in \Phi$. Suppose $0\in \Phi$. Then a $({\cal D}, F,
\Phi)$-approximation is a $\big(\cup_{i\in \Phi}F^i{\cal
D}\big)$-approximation, but the converse is not true in general. If
$0\notin \Phi$, then the two concepts are independent. So, roughly
speaking, a cohomological $\cal D$-approximation with respect to
$(F,\Phi)$ emphasises not only the factorisations but also that the
object belongs to the given subcategory $\cal D$ (and not to
$F^i{\cal D}$ for $0\ne i\in \Phi$).

\section{Derived equivalences for $\Phi$-Yoneda algebras
\label{sect3}}

In this section, we shall prove Theorem \ref{thm1} and derive some
consequences and some simplifications in special cases.

Suppose that $\cal T$ is a triangulated $R$-category with a shift
functor [1], and $M$ is an object in $\cal T$. Suppose that $F$ is
an auto-equivalence of $\cal T$, which may be different from the
shift functor.
%For instance, the derived category of coherent
%sheaves of a projective variety over $\mathbb C$ may have many
%auto-equivalences (see, for example, \cite{st}).

For a subset $\Phi$ of $\mathbb Z$, we define $-\Phi:=\{-x\mid
x\in\Phi\}$, and

$$ {\mathscr X}_{\cal T}^{F,\Phi}(M)=\big\{ X\in {\cal T}\mid
\Hom_{\cal T}(X, F^iM)=0 \mbox{\; for all \;}
i\in\Phi\setminus\{0\}\big\},
$$
$${\mathscr Y}_{\cal T}^{F,\Phi}(M)=\big\{ Y\in {\cal T}\mid \Hom_{\cal T}(M,F^iY)=0 \mbox{\;
for all \;} i\in \Phi\setminus \{0\}\big\}.$$

Let $n$ be a positive integer. For simplicity, we write ${\mathscr
X}^{F,n}(M)$ for ${\mathscr X}_{\cal T}^{F,\{0,1,2,\cdots, n\}}(M)$,
and ${\mathscr X}^{F,\infty}(M)$ for ${\mathscr X}_{\cal
T}^{F,\mathbb N}(M)$ if $\cal T$ is clear in the context. Similarly,
the notations ${\mathscr Y}^{F,n}(M)$ and ${\mathscr
Y}^{F,\infty}(M)$ are defined.

As usual, $F$ is omitted in notation when it coincides with the
shift functor.

Given a triangle $X\stackrel{\alpha}{\lra} M_1\stackrel{\beta}{\lra}
Y \stackrel{w}{\lra} X[1]$ in $\mathcal T$ with $M_1\in \add(M)$, we
define
$$\widetilde{w}= (w,0): Y\lra (X\oplus M)[1], \quad \bar{w}=(0, w)^T: M\oplus Y\lra X[1], $$
 where $(0, w)^T$ stands for the transpose of the matrix
$(0,w)$, and

\medskip $I:=\big\{x=(x_i)\in \Ex{\cal T}{F,\Phi}{X\oplus M}\mid
x_i=0 \; \mbox{for\;} 0\ne i\in \Phi, x_0 \;\mbox{factorises
through}\; \add(M)\;\mbox{and}\; \widetilde{w}[-1] \big\},$

\smallskip $J:=\big\{y=(y_i)\in \Ex{\cal T}{F,\Phi}{M\oplus Y}\mid
y_i=0 \;\mbox{for\;} 0\ne i\in \Phi, y_0 \;\mbox{factorises
through}\; \add(M)\; \mbox{and}\; \bar{w} \big\}.$

\medskip
{\parindent=0pt The} sets $I$ and $J$ are indeed independent of $F$
and $\Phi\setminus \{0\}$, and contained in $\End_{\cal T}(X\oplus
M)$ and $\End_{\cal T}(M\oplus Y)$, respectively.

The main result of this paper is the following theorem which is a
reformulation of Theorem \ref{thm1}.

\begin{Theo} Let $\Phi$ be an admissible subset of $\mathbb Z$,
let $\cal T$ be a triangulated $R$-category with an auto-equivalence
$F$, and let $M$ be an object in $\cal T$. Suppose that
$$X\stackrel{\alpha}{\lra} M_1\stackrel{\beta}{\lra} Y
\stackrel{w}{\lra} X[1]$$ is a triangle in $\mathcal T$ such that
the morphism $\alpha$ is a left $(\add(M), F, \Phi)$-approximation
of $X$, that the morphism $\beta$ is a right $(\add(M), F,
-\Phi)$-approximation of $Y$ and that $X\in {\mathscr
Y}^{F,\Phi}(M)$ and $Y\in {\mathscr X}^{F,\Phi}(M)$. Then the
algebras $\Ex{\cal T}{F,\Phi}{X\oplus M}/I$ and $\Ex{\cal
T}{F,\Phi}{M\oplus Y}/J$ are derived equivalent. \label{verygthm1}
\end{Theo}

{\it Proof.} Let $V= X\oplus M$ and $W=M\oplus Y$. Set
$$\bar{\alpha}:=(\alpha,\,0): X\ra M_1\oplus M, \quad
\bar{\beta}:=\left({{0}\atop{1}}\,\,{{\beta}\atop{0}}\right):
M_1\oplus M\ra M\oplus Y, \quad \bar{w}:=\left({{0}\atop{w}}\right):
M\oplus Y\ra X[1];$$
$$\widetilde{\alpha}:=\left({{\alpha}\atop{0}}\,\,{{0}\atop{1}}\right):
X\oplus M\ra M_1\oplus M,\quad
\widetilde{\beta}:=\left({{\beta}\atop{0}}\right): M_1\oplus M\ra
Y,\quad \widetilde{w}:=(w, 0): Y\lra (X\oplus M)[1].$$ Then there
are two triangles in $\cal T$:
$$\begin{CD}X@>{\bar{\alpha}}>> M_1\oplus
M@>{\bar{\beta}}>> W@>{\bar{w}}>>X[1],
\end{CD} $$
$$\begin{CD}Y[-1]@>{-\widetilde{w}[-1]}>> V@>{\widetilde{\alpha}}>> M_1\oplus M@>{\widetilde{\beta}}>> Y.
\end{CD}$$

Since $F$ is a triangle functor, there is a natural isomorphism
$\delta: F[1]\ra [1]F$. That is, for any object $X$ in $\cal T$,
there is an isomorphism $\delta_X: F(X[1])\ra (FX)[1]$, which is
natural in $X$. The isomorphism $F^i(X[j])\lra (F^iX)[j]$ is denoted
by $\delta(F,i, X, j)$. In part II of this article, there will be
further discussion of these natural transformations.

%First, we have the following lemma.

\begin{Lem} $(1)$ For any morphism $x_i: V\ra F^iV$ with $i\in \Phi$, there is
a morphism $t_i: Y[-1]\ra (F^iY)[-1]$ such that
$(\widetilde{w}[-1])x_i=t_i\delta(F,i,Y,-1)^{-1}\big(F^i(\widetilde{w}[-1])\big).$

$(2)$ For any morphism $y_i: W \ra F^iW$ with $i\in \Phi$, there is
a morphism $t_i: X[1]\ra (F^iX)[1]$ such that
$y_i(F^i\bar{w})\delta(F,i,X,1)$ = $\bar{w}t_i.$ \label{easy}
\end{Lem}

{\it Proof.} (1) Note that $\widetilde{\alpha}$ is a left $(\add(M),
F, \Phi)$-approximation of $V$. Thus, given $x_i: V\ra F^iV$, there
is a morphism $y_i: M_1\oplus M\ra F^i(M_1\oplus M)$ such that
$\widetilde{\alpha}y_i=x_i(F^i\widetilde{\alpha})$. Since $F$ is a
triangle functor, the second triangle implies that there is a
triangle (see \cite[p.4]{HappelTriangle})
$$\begin{CD}
(F^iY)[-1]
@>{\delta(F,i,Y,-1)^{-1}\big(-F^i(\widetilde{w}[-1])\big)}>>F^iV
@>{F^i{\widetilde{\alpha}}}>>F^i(M_1\oplus
M)@>{F^i{\widetilde{\beta}}}>>F^iY.
\end{CD}$$
Thus there is a morphism $t_i: Y[-1]\ra (F^iY)[-1]$ such that
$(\widetilde{w}[-1])x_i=t_i\delta(F,i,Y,-1)^{-1}\big(F^i(\widetilde{w}[-1])\big).$

(2) The proof of (2) is similar to that of (1), using the following
triangle
$$\begin{CD} F^iX@>{F^i\bar{\alpha}}>>F^i(M_1\oplus M)@>{F^i\bar{\beta}}>>F^iW@>{(F^i\bar{w})\delta(F,i,X,1)}>> (F^iX)[1].\end{CD}$$ $\square$

\medskip
Now we prove that the quotient rings in Theorem \ref{verygthm1} are
well-defined.

\begin{Lem} The $I$ and $J$ appearing in {\rm Theorem \ref{verygthm1}} are ideals
of $\Ex{\cal T}{F,\Phi}{V}$ and $\Ex{\cal T}{F,\Phi}{W}$,
respectively. \label{ideal}
\end{Lem}

{\it Proof.} By definition, a morphism $f: V\ra V$ factorises
through $\add(M)$ if and only if there is an object $M'$ in
$\add(M)$ and there are two morphisms $f_1: V\ra M'$ and $f_2: M'\ra
V$ such that $f=f_1f_2$. A morphism $g: V\ra V$ factorises through
$\widetilde{w}[-1]: Y[-1]\ra V$ if and only if there is a morphism
$g': V\ra Y[-1]$ such that $g=g'\big(\widetilde{w}[-1]\big)$. In the
following, we shall prove that $I$ is an ideal in $\Ex{\cal
T}{F,\Phi}{V}$.

The set $I$ is closed under addition in $\Ex{\cal T}{F,\Phi}{V}$. To
show that $I$ is a two-sided ideal in $\Ex{\cal T}{F,\Phi}{V}$, we
pick an $x=(x_i)_{i\in \Phi}\in I$ and a $y=(y_i)_{i\in\Phi}\in
\Ex{\cal T}{F,\Phi}{V}$, and calculate the products $xy$ and $yx$ in
$\Ex{\cal T}{F,\Phi}{V}$. Note that $xy=(x_0y_i)_{i\in \Phi}$ and
$yx=(y_iF^ix_0)_{i\in \Phi}$ since $x_i=0$ for $0\ne i\in \Phi$. We
write $x_0=uv$ for $u: V\ra M'$ and $v: M'\ra V$, where $M'$ is an
object in $\add(M)$, and $x_0= s(\widetilde{w}[-1])$ for a morphism
$s: V\ra Y[-1]$.

We first show that $I$ is a right ideal.

(1) Let $i=0$. The map $x_0y_0$ factorises through an object in
$\add(M)$. Since $x_0$ factorises through $\widetilde{w}[-1]$, it
follows from Lemma \ref{easy} (1) that $x_0y_0$ factorises also
through $\widetilde{w}[-1]$.

(2) Let $0\ne i\in \Phi$. In this case, $\Hom_{\cal T}(M,F^iX)=0$ by
the assumption $X\in {\mathscr Y}^{F,\Phi}(M)$. Let $p_X$ and $p_M$
be the projections of $V$ onto $X$ and $M$, respectively. Then the
composition $vy_iF^ip_X: M'\stackrel{v}{\lra} V\stackrel{y_i}{\lra}
F^iV\stackrel{F^ip_X}{\lra} F^iX$ belongs to $\Hom_{\cal
T}(M',F^iX)=0$. Thus $x_0y_iF^ip_X=uvy_iF^ip_X=0$. By Lemma
\ref{easy} (1), there is a morphism $t_i: Y[-1]\ra F^iY[-1]$ such
that $(\widetilde{w}[-1])y_i=
t_i\delta(F,i,Y,-1)^{-1}F^i(\widetilde{w}[-1])$. Hence
$$\begin{array}{rl} x_0y_i(F^ip_M)& =s (\widetilde{w}[-1])y_i(F^ip_M)=
st_i\delta(F,i,Y,-1)^{-1}F^i(\widetilde{w}[-1])(F^ip_M)\\ & \\ & =
st_i\delta(F,i,Y,-1)^{-1}F^i\big(\widetilde{w}[-1]p_M\big) \\ & \\
&= st_i
\delta(F,i,Y,-1)^{-1}F^i\big((w[-1],0){\tiny\left({{0}\atop{1_M}}\right)}\big)=0.\end{array}$$
Altogether, $x_0y_i=x_0y_i(F^ip_X,F^ip_M)=0$ for $0\ne i\in \Phi$.

Hence $xy\in I$, and $I$ is a right ideal in $\Ex{\cal
T}{F,\Phi}{V}$.

\smallskip
Next, we show that $I$ is a left ideal, that is, we check
$(y_iF^ix_0)_{i\in \Phi}\in I$.

(3) The map $y_0x_0$ factorises through an object in $\add(M)$ and
through $\widetilde{w}[-1]$.

(4) Let $0\ne i\in \Phi$.  Note that $\widetilde{\alpha}: V\ra
M_1\oplus M$ is a left $(\add(M), F, \Phi)$-approximation of $V$.
Thus there is a morphism $h_i:M_1\oplus M\ra F^i(M')$ such that
$y_i(F^iu)={\widetilde{\alpha}h_i}$. By assumption, $\Hom_{\cal
T}(M, F^iX)=0$. This implies that $h_i(F^iv)(F^ip_X)=0$, and
therefore $y_i(F^ix_0)(F^ip_X)=0$. Since
$(F^i\widetilde{w}[-1])(F^ip_M)=0$, we get $y_i(F^ix_0)(F^ip_M)=0$.
Thus $y_iF^ix_0=0$ for $0\ne i\in \Phi$.

Hence $yx\in I$, and $I$ is a left ideal in $\E{\cal T}{F,\Phi}{V}$.
Thus $I$ is an ideal in $\Ex{\cal T}{F,\Phi}{V}$.

Similarly, $J$ is an ideal in $\Ex{\cal T}{F,\Phi}{W}$. $\square$

\medskip
We know that $\E{\cal T}{F,\Phi}{V}{Z}$ is a $\Ex{\cal
T}{F,\Phi}{V}$-module for any object $Z$ in $\cal T$. The next lemma
shows that the ideal $I$ of  $\Ex{\cal T}{F,\Phi}{V}$ may annihilate
some modules of this form.

\medskip
\begin{Lem} Keep the notations as above. Then

\smallskip
$(1)$  $I\cdot \E{\cal T}{F,\Phi}{V}{M}=0$.

\smallskip
$(2)$  $I\cdot \E{\cal T}{F,\Phi}{V}{X}=
\big\{(x_i)_{i\in\Phi}\in\E{\cal T}{F,\Phi}{V}{X}\mid x_i=0 \;
\emph{for\;} 0\ne i\in \Phi, x_0 \;\emph{factorises through}\;
\add(M)\;\emph{and}\; w[-1] \}.$

\smallskip
{ $(3)$} For $x=(x_i)_{i\in\Phi}\in \E{\cal T}{F,\Phi}{V'}{X}$ with
$V'\in \add(V)$, we have $\emph{Im}(\mu(x))\subseteq I\cdot \E{\cal
T}{F,\Phi}{V}{X}$ if and only if $x_i=0$ for all $0\ne i\in \Phi$
and $x_0$ factorises through $\add(M)$ and $w[-1]$, where $\mu$ is
defined in {\rm Lemma \ref{Extprop} (1)}.

$(4)$ Let $f: M'\ra X$ with $M'\in \add(M)$. Then
$\emph{Im}\big(\E{\cal T}{F,\Phi}{V}{f}\big)\subseteq I\cdot\E{\cal
T}{F,\Phi}{V}{X}$ if and only if $f$ factorises through $w[-1]$.
\label{vanish}
\end{Lem}

{\it Proof.} (1) We denote by $\lambda_M=(0,1): M\ra V$ the
canonical inclusion. Let $(x_i)_{i\in \Phi}\in I$ and
$(y_i)_{i\in\Phi}\in \E{\cal T}{F,\Phi}{V}{M}$. Then $(x_i)(y_i)=
(x_0y_i)_{i\in \Phi}$ since $x_i=0$ for $0\ne i\in \Phi$. Since $I$
is an ideal in $\Ex{\cal T}{F,\Phi}{V}$, it follows that
$x\big(y_i(F^i\lambda_M)\big)_{i\in \Phi}=
(x_0y_i(F^i\lambda_M))_{i\in \Phi}\in I$. By the definition of $I$,
we have $x_0y_i(F^i\lambda_M)=0$ for all $0\ne i\in \Phi$ and
$x_0y_0\lambda_M$ factorises through $\widetilde{w}[-1]$. Moreover,
$x_0y_0\lambda_M=(x_0y_0\lambda_Mp_M)\lambda_M
=s\big(\widetilde{w}[-1]p_M\big)\lambda_M= s\cdot 0 \cdot
\lambda_M=0$, where $s$ is a morphism from $V$ to $Y[-1]$. Hence
$x_0y_i(F^i\lambda_M)=0$, and
$x_0y_i=x_0y_i(F^i\lambda_M)(F^ip_M)=0\cdot F^ip_M=0$ for all $i\in
\Phi$. Thus (1) follows.

(2) Let $\lambda_X: X\ra V$ be the canonical inclusion. As in case
(1), it follows that, for $(x_i)_{i\in \Phi}\in I$ and
$(y_i)_{i\in\Phi}\in \E{\cal T}{F,\Phi}{V}{X}$, we have $(x_i)(y_i)=
(x_0y_i)_{i\in \Phi}$, and that $x_0y_0\lambda_X$ factorises through
$\widetilde{w}[-1]$ and $\add(M)$. Hence
$x_0y_0=(x_0y_0\lambda_X)p_X= s(\widetilde{w}[-1])p_X)=s(w[-1])$,
where $s$ is a morphism from $V$ to $Y[-1]$. Conversely, let
$x=(x_i)_{i\in\Phi}\in\E{\cal T}{F, \Phi}{V}{X}$, and suppose that
$x_i=0$ for all $0\neq i\in\Phi$ and that $x_0$ factorises through
$\add(M)$ and $w[-1]$. For $f: U\ra Z$ in ${\cal T}$, we denote by
$\underline{f}$ the element of $\E{\cal T}{F,\Phi}{U}{Z}$
concentrated only in degree $0\in \Phi$. Then it is straightforward
to check that $x\underline{\lambda_X}$ belongs to $I$. Thus,
$x=x\underline{\lambda_X}\underline{p_X}\in I\cdot\E{\cal T}{F,
\Phi}{V}{X}$.

(3) First, suppose $V'=V$ and $\mbox{Im}\big(\mu(x)\big)\subseteq
I\cdot \E{\cal T}{F,\Phi}{V}{X}$. Then $x$, the image of
$\underline{1_V}$ under $\mu(x)$, belongs to $I\cdot \E{\cal
T}{F,\Phi}{V}{X}$. Thus, by (2), we have $x_i=0$ for all $0\ne i\in
\Phi$ and that $x_0$ factorises through $\add(M)$ and $w[-1]$.
Conversely, suppose that $x\in I\cdot\E{\cal T}{F,\Phi}{V}{X}$.
Then, for any $y\in \Ex{\cal T}{F,\Phi}{V}$, the image of $y$ under
$\mu(x)$ is $y\cdot x$. Since $I\cdot\E{\cal T}{F,\Phi}{V}{X}$ is a
$\Ex{\cal T}{F,\Phi}{V}$-submodule of $\E{\cal T}{F,\Phi}{V}{X}$, we
have $yx\in I\cdot\E{\cal T}{F,\Phi}{V}{X}$.

Secondly, suppose that $V'$ is a direct sum of $n$ copies of $V$,
and $x\in \E{\cal T}{F,\Phi}{V'}{X}$. We identify $\E{\cal
T}{F,\Phi}{V'}{X}$ with $\bigoplus_{i=1}^n\E{\cal T}{F,\Phi}{V}{X}$,
and write $x=(a_1, \cdots, a_n)^T$, a column matrix with $a_i\in
\E{\cal T}{F,\Phi}{V}{X}$. Then the image of $\mu(x)$ is the sum of
the image of $\mu(a_i)$ for $1\le i\le n$. Now the conclusion
follows from the first case.

Finally, suppose that $V'$ is a direct summand of $n$ copies of $V$,
that is, $\bigoplus_{i=1}^nV=V'\oplus V''$. If $x\in \E{\cal
T}{F,\Phi}{V'}{X}$, then we may consider $(x,0)^T$ as an element in
$\E{\cal T}{F,\Phi}{\bigoplus_{i=1}^n V}{X}$. Then the proof is
reduced to the second case.

(4) follows from (3) because of $\E{\cal
T}{\Phi}{V}{f}=\mu(\underline{f})$. $\square$

\medskip
Let $\cpx{\widetilde{T}}$ be the complex
$$\begin{CD} \cpx{\widetilde{T}}: \quad 0 @>>> \E{\mathcal T}{F,\Phi}{V}{X}@>{\E{\cal T}{F,\Phi}{V}{\bar{\alpha}}}>>
\E{\cal T}{F,\Phi}{V}{M_1\oplus M}@>>> 0, \end{CD}$$ where the term
$\E{\cal T}{F,\Phi}{V}{X}$ is in degree zero. Then it is the direct
sum of the following two complexes
$$\begin{CD} \quad 0 @>>> \E{\cal T}{F,\Phi}{V}{X}@>{\E{\cal T}{F,\Phi}{V}{\alpha}}>>\E{\cal T}{F,\Phi}{V}{M_1}@>>> 0, \\
0@>>>0 @>>>\E{\cal T}{F,\Phi}{V}{M}@>>> 0.
\end{CD}$$
Let $P=\E{\cal T}{F,\Phi}{V}{X}/I\cdot \E{\cal T}{F,\Phi}{V}{X}$,
and let $p: \E{\cal T}{F,\Phi}{V}{X}\ra P$ be the canonical
surjection. Then, by Lemma \ref{vanish} (1), we may write
$\E{\mathcal T}{F,\Phi}{V}{\bar{\alpha}}=pq$ with $q: P\ra \E{\cal
T}{F,\Phi}{V}{X}$. The complex
$$ \cpx{T}:\qquad
0\lra P\lra \E{\cal}{F,\Phi}{V}{M_1\oplus M}\lra 0$$ in
$\Db{\Ex{\cal T}{F,\Phi}{V}/I}$ is the direct sum of the complexes
$$\begin{CD} 0 @>>> P @>q>>\E{\cal T}{F,\Phi}{V}{M_1}@>>> 0, \\
0@>>>0 @>>>\E{\cal T}{F,\Phi}{V}{M}@>>> 0.
\end{CD}$$
Each term of $\cpx{T}$ is a finitely generated projective $\Ex{\cal
T}{F,\Phi}{V}/I$-module.

\medskip
Before proceeding further, we need to introduce some more notation.
Set
$$ \Lambda := \Ex{\cal T}{F,\Phi}{V},\quad \Gamma := \Ex{\cal
T}{F,\Phi}{W}, \quad \overline{\Lambda} := \Lambda/I, \quad
\overline{\Gamma} := \Gamma/J, $$ where $I$ and $J$ are defined just
before Theorem \ref{verygthm1}.

\begin{Lem} $\cpx{T}$ is a tilting complex over
$\overline{\Lambda}$. \label{tiltcpx}
\end{Lem}

{\it Proof.} It is clear that
$\Hom_{\Kb{\pmodcat{\overline{\Lambda}}}}(\cpx{T},\cpx{T}[i]) =0$
for $i\le -2$ and for $i\ge 2$. We have to check that
$\Hom_{\Kb{\pmodcat{\overline{\Lambda}}}}(\cpx{T},\cpx{T}[1]) =0$
and $\Hom_{\Kb{\pmodcat{\overline{\Lambda}}}}(\cpx{T},\cpx{T}[-1])
=0.$

Let $\cpx{f}\in
\Hom_{\Kb{\pmodcat{\overline{\Lambda}}}}(\cpx{T},\cpx{T}[1])$.
Consider the following diagram:

$$\begin{CD}
@.@. \E{\cal T}{F,\Phi}{V}{X}\\
@.@.@VV{p}V\\
@.0@>>>P@>{q}>>\E{\cal T}{F,\Phi}{V}{M_1\oplus M}@>>>0\\
@.@VVV @VV{f^0}V @VVV\\
0@>>>P@>q>>\E{\cal T}{F,\Phi}{V}{M_1\oplus M}@>>>0
\end{CD}$$
Since both $X$ and $M_1\oplus M$ are in $\add(V)$, Lemma
\ref{Extprop} (1) provides an isomorphism $\mu: \E{\cal
T}{F,\Phi}{X}{M_1\oplus M}\simeq \Hom_{\Lambda}\big(\E{\mathcal
T}{F,\Phi}{V}{X}, \E{\mathcal T}{F,\Phi}{V}{M_1\oplus M}\big)$ and
an element $u=(u_i)_{i\in \Phi}\in\E{\mathcal
T}{F,\Phi}{X}{M_1\oplus M}$ such that $ pf^0 = \mu(u)$. By
assumption, $\bar{\alpha}$ is a left $(\add(M), F,
\Phi)$-approximation of $X$. This yields for each $i\in \Phi$ a
morphism $u_i': M_1\oplus M\ra F^i(M_1\oplus M)$ such that
$u_i=\bar{\alpha}u_i'$. Clearly, $u':=(u_i')_{i\in \Phi}\in
\E{\mathcal T}{F,\Phi}{M_1\oplus M}{M_1\oplus M}$, and $\mu(u')\in
\Hom_{\Lambda}\big(\E{\mathcal T}{F,\Phi}{V}{M_1\oplus M},
\E{\mathcal T}{F,\Phi}{V}{M_1\oplus M}\big)$. Now, we have to check
the following diagram is commutative:
$$ \begin{CD} \E{\mathcal T}{F,\Phi}{V}{X}@>{\E{\mathcal T}{F,\Phi}{V}{\bar{\alpha}}}>> \E{\mathcal T}{F,\Phi}{V}{M_1\oplus M} \\
@V{\mu(u)}VV @VV{\mu(u')}V \\
\E{\mathcal T}{F,\Phi}{V}{M_1\oplus M}@= \E{\mathcal
T}{F,\Phi}{V}{M_1\oplus M}\end{CD}
$$
In fact, if $a=(a_j)_{j\in\Phi}\in \E{\mathcal T}{F,\Phi}{V}{X}$,
then it is sent to $b:=\big(a_jF^j(\bar{\alpha})\big)_{j\in \Phi}$
by $\E{\mathcal T}{F,\Phi}{V}{\bar{\alpha}}$, and further sent to
$bu'=\big(a_j(F^j\bar{\alpha})\big)_{j\in \Phi}u'$ by $\mu(u')$. An
easy calculation shows that $bu'=au$, the image of $a$ under
$\mu(u).$ Thus the diagram is commutative, and
$$ pf^0 = \mu(u) = \E{\mathcal T}{F,\Phi}{V}{\bar{\alpha}}\mu(u')=pq\mu(u').$$
This means that $f^0=q\mu(u')$ - since $p$ is surjective - and that
$\cpx{f}=0$ in ${\Kb{\pmodcat{\overline{\Lambda}}}}$. Therefore
$\Hom_{\Kb{\pmodcat{\overline{\Lambda}}}}(\cpx{T},\cpx{T}[1]) = 0$.

Let $\cpx{f}\in
\Hom_{\Kb{\pmodcat{\overline{\Lambda}}}}(\cpx{T},\cpx{T}[-1])$.
Consider the following diagram:

$$\begin{CD}
0@>>>P@>q>>\E{\cal T}{F,\Phi}{V}{M_1\oplus M}@>>>0\\
@.@VVV @VV{f^1}V @VVV\\
@.0@>>>P@>q>>\E{\cal T}{F,\Phi}{V}{M_1\oplus M}@>>>0
\end{CD}$$
Since $p$ is surjective and $\E{\cal T}{F,\Phi}{V}{M_1\oplus M}$ is
projective in $\Lambda$-mod, $f^1$ can be lifted along $p$, say
$f^1=gp$ with $g: \E{\cal T}{\Phi}{V}{M_1\oplus M}\ra \E{\cal
T}{F,\Phi}{V}{X}$. By assumption, we have $X\in {\mathscr
Y}^{F,\Phi}(M)$, and, by Lemma \ref{Extprop} (3), there is a
homomorphism $u: M_1\oplus M\ra X$ such that $g = \E{\cal
T}{F,\Phi}{V}{u}$. Thus
$$\E{\cal T}{F,\Phi}{V}{u\bar{\alpha}}=\E{\cal T}{F,\Phi}{V}{u}\E{\cal T}{F,\Phi}{V}{\bar{\alpha}}=gpq=f^1q=0. $$
Lemma \ref{Extprop} (2) implies $u\bar{\alpha}=0=u\alpha$. Therefore
$u$ factorises through $-w[-1]$. By Lemma \ref{vanish} (4), the
image of $g \;(=\E{\cal T}{F,\Phi}{V}{u})$ is contained in
$I\cdot\E{\cal T}{F,\Phi}{V}{X}$. It follows that $f^1=gp=0$ and
$\cpx{f}=0$. Hence
$\Hom_{\Kb{\pmodcat{\overline{\Lambda}}}}(\cpx{T},\cpx{T}[-1]) =0.$

Now, $\add(\cpx{T})$ generates $\Kb{\pmodcat{\overline{\Lambda}}}$
as a triangulated category. Thus $\cpx{T}$ is a tilting complex over
$\overline{\Lambda}$.
  $\square$

\medskip
{\it Remark.} To get a tilting complex from $\cpx{\widetilde{T}}$,
one may consider the ideal $I_0$ of $\Ex{\cal T}{\Phi}{V}$
consisting of all endomorphisms $V\ra V$ which are of the form $fg$
with $f: V\ra M'$ and $g: M'\ra V$ such that $M'\in\add(M)$ and
$g\widetilde{\alpha}=0$. Then it is easy to show that the quotient
complex of $\cpx{\widetilde{T}}$ modulo $I_0\cpx{\widetilde{T}}$ is
a two-term tilting complex over $\Ex{\cal T}{\Phi}{V}/I_0$. We will
not use this complex because its endomorphism  algebra cannot be
described in a nice way. Note that the ideal $I_0$ of $\Ex{\cal
T}{\Phi}{V}$ is properly contained in $I$ in general.

\begin{Lem} The two rings $\overline{\Gamma}$ and
$\End_{\Kb{\pmodcat{\overline{\Lambda}}}}(\cpx{T}) $ are isomorphic.
\end{Lem}

{\it Proof.} Since $\overline{\Lambda}$ is a quotient algebra of
$\Lambda$, the category $\modcat{\overline{\Lambda}}$ can be viewed
as a full subcategory of $\modcat{\Lambda}$, and it follows that
$\Kb{\overline{\Lambda}}$ can be viewed as a full subcategory of
$\Kb{\Lambda}$. Thus, we have an isomorphism
$\End_{\Kb{\pmodcat{\overline{\Lambda}}}}(\cpx{T})\simeq
\End_{\Kb{\Lambda}}(\cpx{T})$. To prove the lemma, we shall
construct an isomorphism from $\End_{\Kb{\Lambda}}(\cpx{T})$ to
$\overline{\Gamma}$.

Let $\cpx{f}\in\End_{\Kb{\Lambda}}(\cpx{T})$. Since $p: \E{\cal
T}{F,\Phi}{V}{X}\lra P$ is an epimorphism and $\E{\cal
T}{F,\Phi}{V}{X}$ is a projective $\Lambda$-module, there is a
$\Lambda$-module homomorphism $u^0: \E{\cal T}{F,\Phi}{V}{X}\lra
\E{\cal T}{F,\Phi}{V}{X}$ such that $u^0p=pf^0$. Let $u^1:=f^1$ and
$u^i=0$ for all $i\neq 0,1$. Then it follows from $$u^0\E{\cal
T}{F,\Phi}{V}{\bar{\alpha}}=u^0pq=pf^0q=pqf^1=\E{\cal
T}{F,\Phi}{V}{\bar{\alpha}}u^1$$ that
$\cpx{u}=(u^i)_{i\in\mathbb{Z}}$ is an endomorphism in
$\End_{\Kb{\Lambda}}(\cpx{\widetilde{T}})$. By Lemma \ref{Extprop}
(1), we can assume that  $u^0=\mu(x)$ and $u^1=\mu(y)$ with
$x=(x_i)_{i\in\Phi}\in\Ex{\cal T}{F,\Phi}{X}$ and
$y=(y_i)_{i\in\Phi}\in\Ex{\cal T}{F,\Phi}{M_1\oplus M}$. Now, it
follows from $\E{\cal T}{F,\Phi}{V}{\bar{\alpha}} u^1= u^0\E{\cal
T}{F,\Phi}{V}{\bar{\alpha}}$ that
$$ (\bar{\alpha}y_i)_{i\in
\Phi}=(x_iF^i\bar{\alpha})_{i\in\Phi},\; \mbox{that is,}\;
\bar{\alpha} y_i=x_iF^i\bar{\alpha}\; \mbox{for}\; i\in \Phi.$$ For
each $i\in\Phi$, we can form the following commutative diagram in
${\cal T}$:
$$ (*)\quad\quad\begin{CD} X @>{\bar{\alpha}}>>M_1\oplus M @>{\bar{\beta}}>> W @>{\bar{w}}>> X[1] \\
@V{x_i}VV @VV{y_i}V @VV{h_i}V @VV{x_i[1]}V\\
F^iX @>{F^i\bar{\alpha}}>>F^i(M_1\oplus M) @>{F^i\bar{\beta}}>> F^iW
@>{(F^i\bar{w})\delta(F,i,X,1)}>> (F^iX)[1].\\
\end{CD}  $$
for some morphism $h_i\in \Hom_{\cal T}(W,F^iW)$. Thus, for each
$\cpx{f}\in\End_{\Kb{\Lambda}}(\cpx{T})$, we get an element
$h:=(h_i)_{i\in\Phi}\in\Gamma$ which is $\Ex{\cal T}{F,\Phi}{W}$ by
definition. This leads us to defining the following correspondence:
$$\Theta: \End_{\Kb{{\Lambda}}}(\cpx{T}) \lra
\overline{\Gamma}=\Gamma/J,\qquad \cpx{f}\mapsto h+J.$$

{\em Claim 1.} $\Theta$ is well-defined. \\
{\em Proof.} Suppose that $\cpx{f}\in\End_{\Kb{\Lambda}}(\cpx{T})$
is null-homotopic, that is, there is a map $r: \E{\cal
T}{F,\Phi}{V}{M_1\oplus M}\lra P$ such that $f^0=qr$ and $f^1=rq$.
Since $p$ is surjective and $\E{\cal T}{F,\Phi}{V}{M_1\oplus M}$ is
projective in $\modcat{\Lambda}$, there is a map $s: \E{\cal
T}{F,\Phi}{V}{M_1\oplus M}\lra \E{\cal T}{F,\Phi}{V}{X}$ such that
$sp=r$. Hence $(u^0-pqs)p=u^0p-pqsp=u^0p-pqr=u^0p-pf^0=0$ and
$u^1=rq=spq$. By the assumption $X\in {\mathscr Y}^{F,\Phi}(M)$,
Lemma \ref{Extprop} (3) yields a map $t: M_1\oplus M\lra X$ such
that $s=\E{\cal T}{F,\Phi}{V}{t}=\mu(\underline{t})$. Therefore,
$$\mu(x-\underline{\bar{\alpha}t})p=\big(u^0-\E{\cal T}{\Phi}{V}{\bar{\alpha}}\E{\cal T}{\Phi}{V}{t}\big)p=(u^0-pqs)p=0$$
and $\mu(y-\underline{t\bar{\alpha}})=u^1-spq=0$. Consequently,
${\rm Im} (\mu(x-\underline{\bar{\alpha}t}))\subseteq I\cdot\E{\cal
T}{F,\Phi}{V}{X}$ and $y-\underline{t\bar{\alpha}}=0$. Thus $y_i=0$
for all $0\neq i\in\Phi$ and $y_0=t\bar{\alpha}$. By Lemma
\ref{vanish} (3), we have $x_i=0$ for all $0\neq i\in\Phi$ and
$x_0-\bar{\alpha}t=ab$ for some morphisms $a: X\lra M'$ and $b:
M'\lra X$ with $M'\in\add(M)$. Since $\bar{\alpha}$ is a left
$\add(M)$-approximation of $X$, there is a morphism $c: M_1\oplus
M\lra M'$ such that $a=\bar{\alpha}c$. It follows that
$$x_0=ab+\bar{\alpha}t=\bar{\alpha}cb+\bar{\alpha}t=\bar{\alpha}(cb+t).$$
Now we consider the commutative diagram $(*)$. Suppose $0\neq
i\in\Phi$. Then we have shown that $x_i=y_i=0$. Hence
$\bar{\beta}h_i=y_iF^i\bar{\beta}=0$. This implies that $h_i$
factorises through $\bar{w}$, and, consequently, that $h_i|_{M}=0$
since $\bar{w}|_M=0$. It follows from
$h_i(F^i\bar{w})\delta(F,i,X,1)=\bar{w}(x_i[1])=0$ that $h_i: W\ra
F^iW$ factorises through $F^i(M_1\oplus M)$. Since $Y\in{\mathscr
X}^{F,\Phi}(M)$, we get $h_i|_Y=0$. Altogether, we have shown that
$h_i=0$ for all $0\neq i\in\Phi$. Now consider the diagram $(*)$ in
case $i=0$. First, we have
$\bar{\beta}h_0=y_0\bar{\beta}=t\bar{\alpha}\bar{\beta}=0$, which
means $h_0$ factorises through $\bar{w}$. Second, since
$h_0\bar{w}=\bar{w}(x_0[1])=\bar{w}(\bar{\alpha}[1])(cb+t)[1]=0$,
the morphism $h_0$ factorises through $M_1\oplus M$ which is in
$\add(M)$. Thus, $h\in J$ and $h+J$ is zero in $\overline{\Gamma}$.
This shows that $\Theta$ is well-defined.

{\em Claim 2.} $\Theta$ is injective. \\
{\em Proof.} Suppose that $\Theta(\cpx{f})=h+J=0+J$. Then $h\in J$,
that is, $h_i=0$ for all $0\neq i\in\Phi$, and $h_0$ factorises
through both $\bar{w}$ and $\add(M)$. Suppose $h_0=\bar{w}s$ for a
morphism $s: X[1]\lra W$. For each $0\neq i\in\Phi$, since
$y_iF^i\bar{\beta}=\bar{\beta}h_i=0$, the morphism $y_i: M_1\oplus
M\ra F^i(M_1\oplus M)$ factorises through $F^iX$, and consequently
$y_i=0$ for all $0\neq i\in\Phi$ since
$X\in\mathscr{Y}^{F,\Phi}(M)$. For each $0\neq i\in\Phi$, it follows
from $\bar{w}(x_i[1])=h_i(F^i\bar{w})\delta(F,i,X,1)=0$ that
$x_i[1]$ factorises through $(M_1\oplus M)[1]$, or equivalently, the
morphism $x_i: X\ra F^iX$ factorises through $M_1\oplus M$. Hence
$x_i=0$ for all $0\neq i\in\Phi$ since
$X\in\mathscr{Y}^{F,\Phi}(M)$. Now we consider the case $i=0$.
First, we have
$y_0\bar{\beta}=\bar{\beta}h_0=\bar{\beta}\bar{w}s=0$, which implies
$y_0=t\bar{\alpha}$ for a morphism $t: M_1\oplus M\lra X$. Second,
$(x_0-\bar{\alpha}t)\bar{\alpha}=\bar{\alpha}y_0-\bar{\alpha}t\bar{\alpha}=\bar{\alpha}y_0-\bar{\alpha}y_0=0$.
It follows that $(x_0-\bar{\alpha}t)\alpha=0$, and therefore
$x_0-\bar{\alpha}t$ factorises through $-w[-1]$. Since $h_0: W\ra W$
factorises through $\add(M)$ and since $\bar{\beta}: M_1\oplus M\ra
W$ is a right $\add(M)$-approximation of $W$, we see that $h_0$
factorises through $\bar{\beta}$, say $h_0=r\bar{\beta}$ for some
$r: W\ra M_1\oplus M$. Thus, $\bar{w}(x_0[1]) = h_0\bar{w} =
r\bar{\beta}\bar{w} = 0$, or equivalently, $(-\bar{w}[-1])x_0=0$. It
follows that $x_0$ factorises through $M_1\oplus M$. Since
$\bar{\alpha}t$ also factorises through $M_1\oplus M$, we see that
$x_0-\bar{\alpha}t$ factorises through $\add(M)$. Thus we have shown
that $x_0-\bar{\alpha}t$ factorises through both $\add(M)$ and
$-w[-1]$. Now, by Lemma \ref{vanish} (3), we have ${\rm
Im}\big(\mu(x)-\E{\cal T}{F,\Phi}{V}{\bar{\alpha}t}\big)={\rm
Im}\big(\mu(x-\underline{\bar{\alpha}t})\big)\subseteq I\cdot\E{\cal
T}{F,\Phi}{V}{X}$. Hence
$$p\big(f^0-q\E{\cal T}{F,\Phi}{V}{t}p\big)=u^0p-pq\E{\cal
T}{F,\Phi}{V}{t}p=\big(\mu(x)-\E{\cal
T}{F,\Phi}{V}{\bar{\alpha}t}\big)p=0.$$ This implies that
$f^0=q\big(\E{\cal T}{\Phi}{V}{t}p\big)$ since $p$ is surjective.
Moreover, one can check that $$f^1=u^1=\mu(y)=\E{\cal
T}{F,\Phi}{V}{t}\E{\cal T}{F,\Phi}{V}{\bar{\alpha}}=\big(\E{\cal
T}{F,\Phi}{V}{t}p\big)q.$$ Hence $\cpx{f}$ is null-homotopic, and
consequently $\Theta$ is injective.

{\em Claim 3.} $\Theta$ is surjective. \\
{\em Proof.} Let $h=(h_i)_{i\in \Phi}\in \Gamma$ with $h_i: W\ra
F^iW$ for $i\in \Phi$. Since $\bar{\beta}$ is a right $(\add(M), F,
-\Phi)$-approximation of $W$, we have a morphism $F^{-i}y_i:
F^{-i}(M_1\oplus M)\ra M_1\oplus M$ such that
$\big(F^{-i}\bar{\beta}\big)\big(F^{-i}h_i\big) =
\big(F^{-i}y_i\big)\bar{\beta}$ for $i\in \Phi$. This means that
there is a commutative diagram

$$\begin{CD} X @>{\bar{\alpha}}>>M_1\oplus M @>{\bar{\beta}}>> W @>{\bar{w}}>> X[1] \\
@V{x_i}VV @VV{y_i}V @VV{h_i}V @VV{x_i[1]}V\\
F^iX @>{F^i\bar{\alpha}}>>F^i(M_1\oplus M) @>{F^i\bar{\beta}}>> F^iW
@>{(F^i\bar{w})\delta(F,i,X,1)}>> F^iX[1].\\
\end{CD} $$
Now, define $x:=(x_i)_{i\in \Phi}\in \Ex{\cal T}{F,\Phi}{X}$,
$y:=(y_i)_{i\in\Phi}\in \Ex{\cal T}{F,\Phi}{M_1\oplus M}$;
$u^0:=\mu(x)$, $u^1:=\mu(y)$ and $u^j:=0$ for $j\ne 0,1$. Then
$\cpx{u}:=(u^i)_{i\in {\mathbb Z}}$ belongs to
$\End_{\Kb{\Lambda}}(\cpx{\widetilde{T}})$. Since $u^0: \E{\cal
T}{F,\Phi}{V}{X}\lra \E{\cal T}{F,\Phi}{V}{X}$ takes elements in
$I\cdot\E{\cal T}{F,\Phi}{V}{X}$ to elements in $I\cdot \E{\cal
T}{F,\Phi}{V}{X}$, the image of $I\cdot \E{\cal T}{F,\Phi}{V}{X}$
under the map $u^0p$ is zero, and consequently, there is a unique
map $f^0: P\lra P$ such that $pf^0=u^0p$. Now we have
$$p(f^0q-qu^1)=pf^0q-pqu^1=u^0pq-pqu^1=u^0\E{\cal
T}{F,\Phi}{V}{\bar{\alpha}}-\E{\cal
T}{F,\Phi}{V}{\bar{\alpha}}u^1=0$$ Hence $f^0q=qu^1$ since $p$ is
surjective. Defining $f^1=u^1$ and $f^j=0$ for all $j\neq 0, 1$, we
see that $\cpx{f}=(f^i)_{i\in\mathbb{Z}}$ is an endomorphism in
$\End_{\Kb{\Lambda}}(\cpx{T})$ and $\Theta(\cpx{f})=h+J$. Thus
$\Theta $ is surjective.

{\em Claim 4.} $\Theta$ is an $R$-algebra homomorphism. \\
{\em Proof.} The map $\Theta$ is $R$-linear, so it preserves
addition. For multiplication, we take $\cpx{f}$ and $\cpx{g}$ in
$\End_{\Kb{\Lambda}}(\cpx{T})$. Let $\cpx{u}$ and $\cpx{v}$ be in
$\End_{\Kb{\Lambda}}(\cpx{\widetilde{T}})$ such that $u^0p=pf^0$,
$u^1=f^1$, $v^0p=pg^0$ and $v^1=g^1$. Suppose that $(u^0,
u^1)=\big(\mu(x), \mu(y)\big)$ and $(v^0, v^1)=\big(\mu(x'),
\mu(y')\big)$ with $x, x'\in\Ex{\cal T}{F,\Phi}{X}$ and $y,
y'\in\Ex{\cal T}{\Phi}{M_1\oplus M}$. Let $h:=(h_i)_{i\in\Phi}$ and
$h':=(h'_i)_{i\in\Phi}$ be in $\Gamma$ making the diagram $(*)$
commutative, that is,
$$\genfrac..{0pt}{0}{\bar{\beta}h_i = y_iF^i\bar{\beta}, \quad
\bar{w}(x_i[1])=h_i(F^i\bar{w})\delta(F,i,X,1),}{
\bar{\beta}h'_i=y'_iF^i\bar{\beta}, \quad
\bar{w}(x'_i[1])=h'_i(F^i\bar{w})\delta(F,i,X,1)}$$ for all
$i\in\Phi$. Then, by definition, we have $\Theta(\cpx{f})=h+J$,
$\Theta(\cpx{g})=h'+J$ and
$$\Theta(\cpx{f})\Theta(\cpx{g})=\big(\sum_{\genfrac..{0pt}{3}{i,j\in\Phi}{i+j=k}}h_i(F^ih'_j)\big)_{k\in\Phi}+J.$$
Now we calculate $\Theta(\cpx{f}\cpx{g})$. Let
$\cpx{s}:=\cpx{u}\cpx{v}$. Then $s^0p=pf^0g^0=p(\cpx{f}\cpx{g})^0$,
$s^1=f^1g^1=(\cpx{f}\cpx{g})^1$, and $(s^0,
 s^1)=\big(\mu(xx'), \mu(yy')\big)$, where
$(xx')_k=\displaystyle{\sum_{\genfrac..{0pt}{3}{i,j\in\Phi}{i+j=k}}x_iF^ix'_j}$,
and
$(yy')_k=\displaystyle{\sum_{\genfrac..{0pt}{3}{i,j\in\Phi}{i+j=k}}y_iF^iy'_j}$.
For each $k\in\Phi$, one has to check that
$$(yy')_kF^k\bar{\beta}=\displaystyle{\big(\sum_{\genfrac..{0pt}{3}{i,j\in\Phi}{i+j=k}}y_iF^iy'_j\big)
F^k\bar{\beta}}=\bar{\beta}\displaystyle{\big(\sum_{\genfrac..{0pt}{3}{i,j\in\Phi}{i+j=k}}h_iF^ih'_j\big)}.$$
However, this follows from
$$\begin{array}{rl} y_i(F^iy_j')(F^{i+j}\bar{\beta}) & =
y_iF^i\big(y_j'(F^j\bar{\beta})\big)\\ & \\
 & = y_iF^i(\bar{\beta}h_j')\\ & \\
 & = y_i (F^i\bar{\beta}) (F^ih_j')\\ & \\
 & = \bar{\beta} h_i (F^ih_j').\end{array}. $$
Similarly, for each $k\in\Phi$, we have
$$\big(\displaystyle{\sum_{\genfrac..{0pt}{3}{i,j\in\Phi}{i+j=k}}h_iF^ih'_j}\big)(F^k\bar{w})\delta(F, i, X, 1)=\bar{w}\big((xx')_k[1]\big).$$
This means
$\Theta(\cpx{f}\cpx{g})=\big(\displaystyle{\sum_{\genfrac..{0pt}{3}{i,j\in\Phi}{i+j=k}}h_iF^ih'_j}\big)_{k\in\Phi}+J
=\Theta(\cpx{f})\Theta(\cpx{g})$. Thus $\Theta$ is a ring
homomorphism, and the proof of Theorem \ref{verygthm1} is finished.
$\square$

\medskip Before proceeding, we comment on
the conditions in Theorem \ref{verygthm1}.

({\bf a}) Let $X\stackrel{\alpha}{\lra} M_1\stackrel{\beta}{\lra} Y
\stackrel{w}{\lra} X[1]$ be a triangle in $\cal T$ with $M_1\in
\add(M)$, $X\in {\mathscr Y}^{F,\Phi}(M)$ and $Y\in {\mathscr
X}^{F,\Phi}(M)$. If $\alpha$ is a left $(\add(M), F,
\Phi)$-approximation of $X$, then $\Hom_{\cal T}(X,F^iM)\simeq
\Hom_{\cal T}(M_1, F^iM)$ for $0\ne i\in \Phi$. Similarly, if
$\beta$ is a right $(\add(M), F, -\Phi)$-approximation of $Y$, then
$\Hom_{\cal T}(M, F^iY)=\Hom_{\cal T}(M,F^iM_1)$ for $0\ne i\in
\Phi$. In particular, if $M$ is an $(F,\Phi)$-self-orthogonal object
of $\cal T$, that is, $\Hom_{\cal T}(M,F^iM)=0$ for every $0\ne i\in
\Phi$, and if $\alpha$ is a left $(\add(M), F, \Phi)$-approximation
of $X$ and $\beta$ is a right $(\add(M), F, -\Phi)$-approximation of
$Y$, then $X\in {\mathscr X}^{F,\Phi}(M)$ and $Y\in {\mathscr
Y}^{F,\Phi}(M)$.

\smallskip
({\bf b}) Under the conditions of Theorem \ref{verygthm1}, there are
isomorphisms $\Hom_{\cal T}(X,F^iX)\simeq \Hom_{\cal T}(Y,F^iY)$ for
every $0\ne i\in \Phi$. In fact, this follows from the following
general statement:

Let $\cal T$ be a triangulated category with a shift functor [1].
Suppose that $F$ is a triangle functor from $\cal T$ to itself, and
that $\cal D$ is a full subcategory of $\cal T$. Let $i$ be a
positive integer. Suppose that
$$ X_j\lraf{\alpha_j}D_j\lraf{\beta_j} Y_j\lra X_j[1]$$ is a triangle in $\cal T$, such
that $\alpha_j$ is a left $({\cal D},F, \{i\})$-approximation of
$X_j$, and that $\Hom_{\cal T}(D', F^i(\beta_j)):\Hom_{\cal T}(D',
F^iD_j)\ra \Hom_{\cal T}(D', F^iY_j)$ is surjective for every $D'\in
{\cal D}$ and $j=1,2$. If $\Hom_{\cal T}({\cal
D},F^iX_j)=0=\Hom_{\cal T}(Y_j,F^i{\cal D})$ for $1\le j\le 2$, then
$\Hom_{\cal T}(X_1, F^iX_2)\simeq \Hom_{\cal T}(Y_1,F^iY_2)$.

\smallskip
{\it Proof.} From the given two triangles the following exact
commutative  diagram can be formed:

{\small $$\begin{CD} @.@.\Hom_{\cal T}(D_1, F^iX_2)@>>>\Hom_{\cal
T}(D_1,
F^iD_2) \\
@.@. @VVV @V{(\alpha_1,F^iD_2)}VV\\
 @.@.\Hom_{\cal T}(X_1, F^iX_2)@>>>\Hom_{\cal T}(X_1,
F^iD_2) \\
@.@. @VVV @V0VV\\
\Hom_{\cal T}(Y_1, F^iD_2) @>>>\Hom_{\cal T}(Y_1,
F^iY_2)@>>>\Hom_{\cal T}(Y_1, F^iX_2[1])@>>>\Hom_{\cal T}(Y_1,
F^iD_2[1]) \\
@VVV @VVV @VV{\qquad\qquad \quad(*)}V @VVV\\
\Hom_{\cal T}(D_1, F^iD_2) @>{(D_1,F^i(\beta_2))}>>\Hom_{\cal
T}(D_1, F^iY_2)@>0>>\Hom_{\cal T}(D_1, F^iX_2[1])@>>>\Hom_{\cal
T}(D_1,
F^iD_2[1]). \\
\end{CD}$$}
Since $\Hom_{\cal T}(Y_1, F^iD_2)=\Hom_{\cal T}(D_1, F^iX_2)=0$ by
assumption and since $\Hom_{\cal T}(\alpha_1,F^iD_2)$ and
$\Hom_{\cal T}(D_1,F^i\beta_2)$ are surjective by the property of
approximation, the conclusion follows from the commutative square
($*$). $\square$

\medskip
({\bf c}) Let $X\stackrel{\alpha}{\lra} M_1\stackrel{\beta}{\lra}
Y\lraf{w} X[1]$ be an $\add(M)$-split triangle in $\cal T$. Define
$V:=X\oplus M$, $\Lambda_0:=\End_{\cal
 T}(V)$, $W:=M\oplus Y$, and $\Gamma_0:=\End_{\cal T}(W)$. Let $I$
 and $J$ be as defined in {\rm Theorem \ref{verygthm1}}. Then the ideals
$I$ and $J$ in Theorem \ref{thm1} have the following
characterisation:

(i) Let $e$ be the idempotent in $\Gamma_0$ corresponding to the
direct summand $M$ of $W$. Then $J$ is the submodule of the left
$\Gamma_0$-module $\Gamma_0e\Gamma_0$, which is maximal with respect
to $eJ=0$.

(ii) Let $f$ be the idempotent in $\Lambda_0$ corresponding to the
direct summand $M$ of $V$. Then $I$ is the submodule of the right
$\Lambda_0$-module $\Lambda_0f\Lambda_0$ which is maximal with
respect to $If=0$.

\smallskip
{\it Proof.} By Lemma \ref{ideal}, the sets $I$ and $J$ are ideals
of $\Lambda_0$ and $\Gamma_0$, respectively.

(i) Let $p_M: W\ra M$ and $\lambda_M: M\ra W$ be the canonical
projection and injection, respectively. By definition,
$e=p_M\lambda_M$. The set $\Gamma_0e\Gamma_0$ is precisely the set
of all endomorphisms of $W$ that factorise through $\add(M)$. The
endomorphisms of $W$ factorising through $\bar{w}$ are those
endomorphisms $x$ that satisfy $\bar{\beta}x=0$, and consequently
$ex=p_M\lambda_Mx=p_M(\bar{\beta}|_{M})x=0$. Hence $J$ is a
submodule of $_{\Gamma_0}\Gamma_0e\Gamma_0$ with $eJ=0$. Suppose
that  $\bar{J}\subseteq {}_{\Gamma_0}\Gamma_0e\Gamma_0$ is another
submodule containing $J$ with $e\bar{J}=0$. Then $e\bar{J}=0$
implies $\Hom_{\Gamma_0}(\Hom_{\cal T}(W, M), \bar{J})=0$, and
consequently $\Hom_{\Gamma_0}(\Hom_{\cal T}(W, M'), \bar{J})=0$ for
all $M'\in\add(M)$. For each $x\in\bar{J}$, the image of the
morphism $\Hom_{\cal T}(W, x)$ is contained in $\bar{J}$ since
$\bar{J}$ is a left ideal of $\Gamma_0$. Thus, the morphism
$\Hom_{\cal T}(W, \bar{\beta}x)$ is a $\Gamma_0$-module morphism
from $\Hom_{\cal T}(W, M_1\oplus M)$ to the image of  $\Hom_{\cal
T}(W, x)$. Hence $\Hom_{\cal T}(W, \bar{\beta}x)=0$, and
consequently $\bar{\beta}x=0$. This implies $x\in J$.  This proves
(i).

(ii) The proof is similar to that of (i). $\square$

\medskip
A special case of Theorem \ref{verygthm1} is the following
corollary.

\begin{Koro}
Let $\Phi$ be an admissible subset of $\mathbb Z$, and let $\cal T$
be a triangulated $R$-category with an auto-equivalence $F$, and let
$M$ be an object in $\cal T$. Suppose that $X\stackrel{\alpha}{\lra}
M_1\stackrel{\beta}{\lra} Y \stackrel{w}{\lra} X[1]$ is an
$\add(M)$-split triangle in $\mathcal T$, and suppose that $X$ and
$Y$ both are in ${\mathscr X}^{F,\Phi}(M)\cap{\mathscr
Y}^{F,\Phi}(M)$. Then $\Ex{\cal T}{F,\Phi}{X\oplus M}/I$ and
$\Ex{\cal T}{F,\Phi}{M\oplus Y}/J$ are derived equivalent.
\label{cortoverygthm1}
\end{Koro}

The following special case of Theorem \ref{verygthm1} is useful to
construct explicit examples.

\begin{Koro}
Let $\cal T$ be a triangulated $R$-category with $[1]$ the shift
functor, and let $M$ be an object in $\cal T$. Suppose that
$X\stackrel{\alpha}{\lra} M_1\stackrel{\beta}{\lra} Y
\stackrel{w}{\lra} X[1]$ is a triangle in $\mathcal T$ such that
$M_1\in \add(M)$, and suppose that $X\in {\mathscr Y}^{n+1}(M)$ and
$Y\in {\mathscr X}^{n+1}(M)$. Then, for any admissible subset $\Phi$
of ${\mathbb N}_n$, the algebras $\Ex{\cal T}{\Phi}{X\oplus M}/I$
and $\Ex{\cal T}{\Phi}{M\oplus Y}/J$ are derived equivalent.
\label{3.8}
\end{Koro}

{\it Proof.} We show that $\beta$ is a right $(\add(M),
-\Phi)$-approximation of $Y$. Note that, for $i\in \Phi$, we always
have $i+1\le n+1$. Hence $\Hom_{\cal T}(M,X[i+1])=0$ for $i\in
\Phi$. Now apply $\Hom_{\cal T}(M[-i],-)$ with $i\in \Phi$ to the
triangle $X\stackrel{\alpha}{\lra} M_1\stackrel{\beta}{\lra} Y
\stackrel{w}{\lra} X[1]$:
$$ \cdots\ra \Hom_{\cal T}(M[-i],M_1)\lra \Hom_{\cal T}(M[-i],Y)\lra\Hom_{\cal T}(M[-i], X[1])\ra\cdots $$
Because of $\Hom_{\cal T}(M[-i],X[1])= \Hom_{\cal T}(M, X[i+1])=0$,
the map $\beta$ is a right $(\add(M), -\Phi)$-approximation of $Y$.

Similarly, it follows from $\Ext^{i+1}_{\cal T}(Y,M)=0$ for $i\in
\Phi$ that $\alpha$ is a left $(\add(M),\Phi)$-approximation of $X$.
Now Corollary \ref{3.8} follows from Theorem \ref{verygthm1}.
$\square$

\medskip
Another special case of Theorem \ref{verygthm1} is that $I=0$ and
$J=0$. Here is a condition when the ideals $I$ and $J$ in Theorem
\ref{verygthm1} vanish.

\begin{Prop} Let $X\stackrel{\alpha}{\lra} M_1\stackrel{\beta}{\lra} Y\lraf{w} X[1]$ be an $\add(M)$-split
triangle in $\cal T$. Define $V:=X\oplus M$, $\Lambda_0:=\End_{\cal
 T}(V)$, $W:=M\oplus Y$, and $\Gamma_0:=\End_{\cal T}(W)$. Let $I'$ be the ideal of $\Lambda_0$ consisting of all
$f: V\ra V$ that factorises through $\widetilde{w}[-1]: Y[-1]\ra V$,
and let $J'$ be the ideal of $\Gamma_0$ consisting of all $g: W\ra
W$ that factorises through $\bar{w}: W\ra X[1]$.

$(1)$  Suppose that $\Lambda_0$ is an Artin algebra.  If
$\add\big(\top_{\Lambda_0}\Hom_{\cal T}(V, X)\big)\cap
\add\big(\top(_{\Lambda_0}D\Lambda_0)\big)=0$, then $I'=0$.

$(2)$ Suppose that $\Gamma_0$ is an Artin algebra.  If
$\add\big(\top_{\Gamma_0}\Hom_{\cal T}(W,Y)\big)\cap
\add\big(\soc({}_{\Gamma_0}{\Gamma_0})\big)=0$, then $J'=0$.
\label{bigideal}
\end{Prop}

{\it Proof.} We prove (1).  The proof of (2) is similar to that of
(1), and we omit it.

We have a triangle
$Y[-1]\lraf{-\widetilde{w}[-1]}V\lraf{\widetilde{\alpha}}M_1\oplus
M\lraf{\widetilde{\beta}} Y$, apply $\Hom_{\cal T}(-,V)$ to this
triangle, and get the following exact sequence of right
$\Lambda_0$-modules:
$$\Hom_{\cal T}(M_1\oplus M,V)\lra \Hom_{\cal T}(V,V)\lra C\lra 0,$$
where $C$ is the cokernel of $\Hom_{\cal T}(\widetilde{\alpha}, V).$
Now, applying $\Hom_{\Lambda_0\opp}(\Hom_{\cal T}(M,V),-)$ to the
above exact sequence, we get another exact sequence which is
isomorphic to the following exact sequence:
$$ \Hom_{\cal T}(M_1\oplus M, M)\lraf{(\widetilde{\alpha},M)}
\Hom_{\cal T}(V,M)\lra \Hom_{\Lambda_0\opp}\big(\Hom_{\cal
T}(M,V),C\big)\lra 0. $$ Since $\widetilde{\alpha}$ is a left
$\add(M)$-approximation of $V$, the map $\Hom_{\cal
T}(\widetilde{\alpha},M)$ is surjective, and consequently
$\Hom_{\Lambda_0\opp}\big(\Hom_{\cal T}(M,V),C\big)=0.$ So, the
right $\Lambda_0$-module $C$ has no composition factors in
$\top\big(\Hom_{\cal T}(M,V)\big)$, and that $C$ has composition
factors only in $\top\big(\Hom_{\cal T}(X,V)\big)$.  This is
equivalent to saying that the $\Lambda_0$-module $D(C)$ has
composition factors only in $\soc\big(D\Hom_{\cal T}(X,V)\big)$
which is isomorphic to $\top\big(\Hom_{\cal T}(V,X)\big)$.

Let $x: V\ra V$ be an element in $I'\subseteq \Lambda_0$. Then $x$
factorises through $-\widetilde{w}[-1]$, or equivalently,
$x\widetilde{\alpha}=0$. This implies that $\big(D\Hom_{\cal
T}(x,V)\big)\big(D\Hom_{\cal T}(\widetilde{\alpha},V)\big)=0$. Thus
the image of $D\Hom_{\cal T}(x,V)$ is contained in the kernel of
$D\Hom_{\cal T}(\widetilde{\alpha},V)$, which is isomorphic to
$D(C)$. Therefore, if $D\Hom_{\cal T}(x,V)\ne 0$, then the top of
the image of $D\Hom_{\cal T}(x,V)$ is contained in
$\add\big(\top_{\Lambda_0}\Hom_{\cal T}(V, X)\big)\cap
\add\big(\top(_{\Lambda_0}D\Lambda_0)\big)=0$, this is a
contradiction. Thus we must have $\Hom_{\cal T}(x,V)=0$. Since
$\Hom_{\cal T}(-,V)$ is a duality from $\add(V)$ to
$\pmodcat{\Lambda_0\opp}$, we obtain $x=0$. Thus $I'=0$. $\square$

\medskip
\emph{Remark}. (1) if we substitute ``$\add(M)$-split" for ``left $
(\add(M),\Phi)$-approximation" and ``right
$(\add(M)$, $-\Phi)$-approximation" in Proposition \ref{bigideal}, and
if we consider $\Ex{\cal T}{\Phi}{V}$ and $\Ex{\cal T}{\Phi}{W}$
instead of $\Lambda_0$ and $\Gamma_0$, then Proposition
\ref{bigideal} is still true. The proof is almost the same.

(2) By definition, there are inclusions $I\subseteq I'$ and
$J\subseteq J'$. Sometimes it is easy to verify that $I'$ and $J'$
vanish if the algebras $\Lambda_0$ and $\Gamma_0$ are described by
quivers with relations.

\medskip
For the derived category of an abelian category, the following
result provides an explicit example for $I = 0 = J$.

\begin{Prop} Let $\cal A$ be an abelian category, and let $M$ be an
object of $\cal A$. Suppose that $0\ra X\stackrel{\alpha}{\lra}
M_1\stackrel{\beta}{\lra} Y\ra 0$ is an exact sequence in $\cal A$
with $M_1\in \add(M)$. Consider the induced triangle
$X\stackrel{\alpha}{\lra} M_1\stackrel{\beta}{\lra} Y\lraf{w} X[1]$
in $\Db{\cal A}$. Then the ideals $I$ and $J$ defined in {\rm
Theorem \ref{verygthm1}} vanish. \label{vanish2}
\end{Prop}

{\it Proof.} Every exact sequence $0\ra X\ra M_1\ra Y\ra 0$ in $\cal
A$ gives rise to a triangle $ X\ra M_1\ra Y\ra X[1]$ in $\Db{\cal
A}$. Now we show that the exactness of the given sequence in $\cal
A$ implies that the two ideals $I$ and $J$ in Theorem
\ref{verygthm1} are equal to zero. Since $I$ is contained in
$\End_{\Db{\cal A}}(X\oplus M)$, it is sufficient to show that if a
morphism $x: X\oplus M\ra X\oplus M$ factorises through $\add(M)$
and $\widetilde{w}[-1]$, then $x=0$. Let $x$ be such a morphism.
Then we see immediately that $x\widetilde{\alpha}=0$ in $\Db{\cal
A}$. Since $\cal A$ is fully embedded in $\Db{\cal A}$, we also have
$x\widetilde{\alpha}=0$ in $\cal A$. Consequently, $x=0$ since
$\widetilde{\alpha}$ is injective in $\cal A$. Thus $I=0$. Dually,
we can show $J=0$. Hence Proposition \ref{vanish2} holds true.
$\square$

\medskip
As an immediate application of the proof of Theorem \ref{verygthm1}
together with a result on derived equivalences in \cite{PX}, we have
the following corollary.

\begin{Koro} We keep all assumptions of {\rm Theorem \ref{verygthm1}}.
If $\overline{\Lambda}$ and $\overline{\Gamma}$ both are left
coherent rings $($for example, if $\Phi$ is finite and $\cal T$=
$\Db{A}$ with $A$ a finite dimensional algebra over a field$)$, then
$\findim(\overline{\Lambda})-1\le \findim(\overline{\Gamma})\le
\findim(\overline{\Lambda})-1$, where $\findim(\overline{\Lambda})$
stands for the finitistic dimension of $\overline{\Lambda}$.
\end{Koro}

Recall that, given a ring $S$ with identity, the \emph{finitistic
dimension} of $S$ is defined to be the supremum of the projective
dimensions of finitely generated $S$-modules of finite projective
dimension.

\medskip
Since the map $q$ in the proof of Theorem \ref{verygthm1} is not
always injective, the tilting complex $\cpx{T}$ is not, in general,
isomorphic in $\Db{\Ex{\cal T}{F,\Phi}{V}/I}$ to a tilting module.
Thus the derived equivalence presented in Theorem \ref{verygthm1} is
not given by a tilting module in general (in contrast with the
situation of Theorem \ref{deadss}). In fact, it is easy to see that
the derived equivalence in Theorem \ref{verygthm1} is given by a
tilting module if the kernel of $\E{\cal T}{F,\Phi}{V}{\alpha}$ is
$I\cdot\E{\cal T}{F,\Phi}{V}{X}$.

Moreover, a small additive category may be embedded into an abelian
category of coherent functors (see \cite[Chapter IV, Section
2]{mitchell}). This will, however, not in general turn a $\cal
D$-split sequence in the additive category into an exact sequence in
the abelian category since otherwise the sequence would split, and
therefore cannot provide a triangle in the derived category of the
abelian category. Consequently, Theorem \ref{deadss} cannot be
obtained from Theorem \ref{verygthm1} by taking $\Phi=\{0\}$ and
embedding an additive category into an abelian category.

Finally, we mention that Theorem \ref{verygthm1} generalises the
result \cite[Proposition 5.1]{hx2} by choosing $\Phi=\{0\}$. Indeed,
under the conditions of \cite[Proposition 5.1]{hx2}, the ideals $I$
and $J$ in Theorem \ref{verygthm1} vanish. Theorem \ref{verygthm1}
covers various other situations, some of which will be discussed in
the next
section. \\
%In the forthcoming second part of this article, variations of
%Theorem \ref{verygthm1} will be given that cover more general and
%more abstract setups.

\section{$\Phi$-Yoneda algebras in some explicit situations
 \label{sect4}}
In this section, we shall describe some natural habitats for Theorem
\ref{verygthm1} and relate it to several widely used concepts that
fit with or simplify the assumptions of Theorem \ref{verygthm1}. Throughout,
we choose $F$ to be the shift functor of the triangulated category
considered.

We note that Alex Dugas, in independent work \cite{Dugas} that also
is motivated by \cite{hx2}, has constructed derived equivalent pairs
of symmetric algebras. As explained in \cite{Dugas} (Remark (3) in
section 4) his examples appear in our framework, too.

\subsection{Derived categories of Artin algebras}

A first consequence of Theorem \ref{verygthm1} is the following
result for ${\cal T}=\Db{A}$ with $A$ an Artin $R$-algebra.

\begin{Theo}
Let $\Phi$ be an admissible subset of $\mathbb N$, let $M$ be an
$A$-module, and let $0\ra X\stackrel{\alpha}{\lra}
M_1\stackrel{\beta}{\lra} Y \ra 0$ be an exact sequence in
$\modcat{A}$ with $\alpha$ a left $(\add(M), \Phi)$-approximation
 of $X$ and $\beta$ a right $(\add(M), -\Phi)$-approximation
 of $Y$ in $\Db{A}$ such that $X\in {\mathscr
Y}^{\Phi}(M)$ and $Y\in {\mathscr X}^{\Phi}(M)$.  Then the
perforated Yoneda algebras $\Ex{A}{\Phi}{X\oplus M}$ and
$\Ex{A}{\Phi}{M\oplus Y}$ are derived equivalent. \label{gthm1}
\end{Theo}

{\it Proof.} This is a consequence of Theorem \ref{verygthm1} and
Proposition \ref{vanish2} if we take $\cal T$ = $\Db{A}$. $\square$

\medskip
Under the assumptions of Theorem \ref{gthm1}, the higher cohomology
groups $\Ext^i_A(X,X)$ of $X$ is isomorphic to the higher cohomology
groups $\Ext_A^i(Y,Y)$ of $Y$ for each $0\ne i\in \Phi$. This
follows from the comment (b) before Corollary \ref{cortoverygthm1}.

When requiring additional orthogonality conditions on $X$ and $Y$ in
Theorem \ref{gthm1}, we get the following corollary.

\begin{Koro}
Let $\Phi$ be an admissible subset of $\mathbb N$, let $M$ be an
$A$-module, and let $0\ra X\stackrel{\alpha}{\lra}
M_1\stackrel{\beta}{\lra} Y \ra 0$ be an  $\add(M)$-split sequence
in $\modcat{A}$ such that $X, Y \in {\mathscr
X}^{\Phi}(M)\cap{\mathscr Y}^{\Phi}(M)$.  Then the perforated Yoneda
algebras $\Ex{A}{\Phi}{X\oplus M}$ and $\Ex{A}{\Phi}{M\oplus Y}$ are
derived equivalent. \label{thm1adss}
\end{Koro}

{\it Proof.} This follows immediately from Corollary
\ref{cortoverygthm1} and Proposition \ref{vanish2}.
 $\square$

\medskip
If the orthogonality conditions in Corollary \ref{thm1adss} hold for
${\mathbb N}_n$ or $\mathbb N$, then we get the following
consequence.

\begin{Koro} Suppose that $M$ is an
$A$-module. Let $0\ra X\stackrel{\alpha}{\lra}
M_1\stackrel{\beta}{\lra} Y \ra 0$ be an  $\add(M)$-split sequence
in $\modcat{A}$ such that $X, Y \in {\mathscr X}^n(M)\cap{\mathscr
Y}^n(M)$ for $n$ a positive number or infinity. Then, for any
admissible subset $\Phi$ of ${\mathbb N}_n$, the perforated Yoneda
algebras $\Ex{A}{\Phi}{X\oplus M}$ and $\Ex{A}{\Phi}{M\oplus Y}$ are
derived equivalent. \label{corthm1adss}
\end{Koro}

\medskip
The following result shows that the orthogonality conditions are
related to the concepts of {\em short cycle} and {\em short chain}
in $A$-mod (see \cite[Chapter IX, p.313]{arsbook}). Recall that a
short cycle of length $2$ from an indecomposable module $X$ to $X$
is a sequence of non-zero radical homomorphisms
$X\lraf{f}M\lraf{g}X$ with $M$ indecomposable; and a short chain is
a sequence of non-zero radical homomorphisms
$X\lraf{f}M\lraf{g}D\Tr(X)$ with $X$ indecomposable .

\begin{Koro} Let
$A$ be an Artin algebra, and let $0\ra X\ra M\ra Y\ra 0$ be an
Auslander-Reiten sequence in $\modcat{A}$. Suppose neither $X$ nor
$Y$ lies on a short cycle of length $2$, nor on a short chain. Then
the trivial extension of $\End_A(X\oplus M)$ by the bimodule
$\Ext^1_A(X, X)\oplus \Ext^1_A(M,M)$ is derived equivalent to the
trivial extension of $\End_A(M\oplus Y)$ by the bimodule
$\Ext^1_A(Y, Y)\oplus \Ext^1_A(M,M)$. \label{shortcycle}
\end{Koro}

{\it Proof.} An Auslander-Reiten sequence $0\ra X\ra M\ra Y\ra 0$ is
always an  $\add(M)$-split sequence. Since $Y$ does not lie on a
short cycle, the Auslander-Reiten formula $D\underline{\Hom}_A(\Tr
D(X), M) \simeq \Ext^1_A(M,X)\simeq D\overline{\Hom}_A(X,D\Tr(M))$
(see \cite[p.131]{arsbook}) implies $\Ext^1_A(M,X)=0$. Moreover, $X$
not lying on a short cycle implies $\Ext^1_A(Y,M)=0$. Similarly, the
Auslander-Reiten formula yields that $\Ext^1_A(X,M)=0$ - since $X$
does not lie on a short chain - and that $\Ext^1_A(M,Y)=0$ - since
$Y$ does not lie on a short chain. Thus Corollary \ref{shortcycle}
follows from Corollary \ref{corthm1adss} when $n=1$. $\square$

\medskip
The next corollary is a consequence of Corollary \ref{corthm1adss}.

\begin{Koro} Let
$A$ be an Artin algebra, and let $X$ be an $A$-module such that
$\Ext^i_A(X,A)=0$ for all $1\le i< n+2$ with $n$ a fixed positive
integer or infinity.  Then, for any admissible subset $\Phi$ of
${\mathbb N}_n$, the perforated Yoneda algebras
$\Ex{A}{\Phi}{A\oplus X}$ and $\Ex{A}{\Phi}{A\oplus \Omega(X)}$ are
derived equivalent. \label{gnc}
\end{Koro}

{\it Proof.} If $\Ext^i_A(X,A)=0$ for a fixed $i\ge 1$, then $0\ra
\Omega^i(X)\ra P_{i-1}\ra \Omega^{i-1}(X)\ra 0$ is  an
$\add(_AA)$-split sequence in $A$-mod, where $P_i$ is a projective
cover of $\Omega^i(X)$. Using this fact, Corollary \ref{gnc} follows
immediately from Corollary \ref{corthm1adss}. $\square$

\medskip
The condition $\Ext^i_A(X,A)=0$ on $X$ in Corollary \ref{gnc} is
related to the context of the {\em Generalised Nakayama Conjecture}.
This states that if an $A$-module $T$ satisfies $\Ext^i_A(A\oplus T,
A\oplus T)=0$ for all $i>0$ then $T$ should be projective. The above
Corollary \ref{gnc} (or \cite[Theorem 1.1]{hx2}) describes the shape
of the syzygy modules $\Omega^i(X)$: If $X$ is indecomposable and
non-projective and satisfies $\Ext^i_A(X,A)=0$ for all $i>0$, then,
for each $j\ge 0$, there is an indecomposable non-projective module
$L_j$ such that $\Omega^j(X)\simeq L_j^{m_j}$ for an integer
$m_j>0$.

In Corollary \ref{gnc}, there are isomorphisms $\Ext^i_A(X,X)\simeq
\Ext^i_A(\Omega(X),\Omega(X))$ for all $i\ge 1$. Thus the algebras
$\Ex{A}{\Phi}{A\oplus X}$ and $\Ex{A}{\Phi}{A\oplus \Omega(X)}$ are
the extensions of $\End_A(A\oplus X)$ and $\End_A(A\oplus
\Omega(X))$ by the same ideal $\E{A}{\Phi\setminus \{0\}}{X}{X}$,
respectively. The algebras $\Ex{A}{\Phi}{X\oplus M}$ and
$\Ex{A}{\Phi}{M\oplus Y}$ in Corollary \ref{corthm1adss}, however,
are the extensions of $\End_A(X\oplus M)$ and $\End_A(M\oplus Y)$ by
possibly different ideals $\Ex{A}{\Phi\setminus \{0\}}{M}\oplus
\Ex{A}{\Phi\setminus \{0\}}{X}$ and $\Ex{A}{\Phi\setminus
\{0\}}{M}\oplus \Ex{A}{\Phi\setminus \{0\}}{Y}$, respectively.

\medskip
Recall that a module $M\in A$-mod is called \emph{reflexive} if the
evaluation map $$\alpha_M: M \ra M^{**} :=
\Hom_{A\opp}\big(\Hom_A(M,A),A_A\big)$$ is an isomorphism of
modules.

\begin{Koro} Let $M$ be a reflexive $A$-module. Then, for any subset
$0\in \Phi\subseteq \{0,1\}$, the perforated Yoneda algebras
$\Ex{A}{\Phi}{D(A_A)\oplus D\Tr(M)}$ and $\Ex{A}{\Phi}{D(A_A)\oplus
\Omega^{-1}(D\Tr(M))}$ are derived equivalent, where $\Omega^{-1}$
is the co-syzygy operator. \label{reflexive}
\end{Koro}

{\it Proof.} By \cite[IV, Proposition 3.2]{arsbook}, the kernel and
cokernel of the evaluation map $\alpha_M$ are
$\Ext^1_{A\opp}(\Tr(M),A)$ and $\Ext^2_{A\opp}(\Tr(M),A)$,
respectively. As $\Ex{A}{\Phi}{U}\simeq \Ex{A\opp}{\Phi}{D(U)}\opp$
for any $A$-module $U$, Corollary \ref{reflexive} follows from
Corollary \ref{gnc} for right modules. $\square$

%\medskip
%Note that the operator $\Omega\Tr$ was introduced  and studied in
%\cite{AB} and denoted by $\lambda$ in \cite{ms} where it was used to
%define the notion of linkage of modules.

\medskip
A special case of Corollary \ref{gnc}, is a result on self-injective
algebras that has been obtained in \cite[Corollary 3.14]{HX4}):

\begin{Koro} If $A$ is a
self-injective Artin algebra, then, for any admissible subset $\Phi$
of ${\mathbb N}$, the perforated Yoneda algebras
$\Ex{A}{\Phi}{A\oplus X}$ and $\Ex{A}{\Phi}{A\oplus \Omega(X)}$ are
derived equivalent.\label{selfinjalg}
\end{Koro}

Another concept related to the Generalised Nakayama Conjectures and
to modules being projective and injective, is the {\em dominant
dimension} of an algebra or a module.

Suppose that $A$ is an Artin $R$-algebra. By definition, the {\em
dominant dimension} of $A$ is greater than or equal to $n$ if in the
minimal injective resolution of $_{A}A$:
$$0\lra A\lra I_0\lra I_1\lra \cdots \lra I_{n-1}\lra I_n \lra \cdots, $$
the first $n$ injective $A$-modules $I_0$, $\cdots, I_{n-1}$ are
projective. In this case we write $\domdim(A)\ge n$. Let $C_i$ be
the $i$-th cosyzygy of $A$, that is, the cokernel of the map
$I_{i-1}\ra I_i$.

For an $A$-module $X$, we define $a(X)$ to be the number of
non-isomorphic indecomposable direct summands of $M$. The
\emph{self-injective measure} of $A$ is defined to be the number
$m(A):= a(A)-a(I_0)$, where $I_0$ is an injective hull of $A$. Thus,
if $A$ is self-injective, then $m(A)=0$. If $\domdim{(A)}\ge 1$,
then $A$ is self-injective if and only if $m(A)=0$. So the Nakayama
conjecture can be reformulated as: If $\domdim{(A)}=\infty$, then
$m(A)=0$.

\begin{Koro} Let $A$ be an Artin algebra, and let $T$ be the direct sum of all non-isomorphic indecomposable
projective-injective $A$-modules.

$(1)$ If $\domdim{(A)}\ge n\ge 2$, then $\End_{A}(T\oplus C_i)$ is
derived equivalent to $A$ for $1\le i< n$.

$(2)$ If $\domdim{(A)}\ge n+1<\infty$, then $m(A)= a(C_n)$.
\label{dom}
\end{Koro}

{\it Proof.} Since  the sequence $0\ra C_{i-1}\ra I_i\ra C_i\ra 0$
is an  $\add(I_i)$-split sequence (or an  $\add(T)$-split sequence),
the orthogonality conditions in Corollary \ref{corthm1adss} are
trivially satisfied.  Derived equivalence preserves the number of
non-isomorphic simple modules. Therefore, Corollary \ref{dom}
follows now from Corollary \ref{corthm1adss}. Here we also use the
observation that $\add(C_i)\cap \add(I_j) =\{0\}$ for all $0\le
i,j\le n$.  Alternatively, one can also use Lemma \ref{deadss} to
prove this corollary. $\square$

\medskip
Examples of algebras of dominant dimension at least $n$ can be
obtained in the following way: Let $A$ be a self-injective algebra
and $X$ an $A$-module. If $\Ext^i_A(X,X)=0$ for all $1\le i\le n$,
then dom.dim$\big(\End_A(A\oplus X)\big)\ge n+2$.
\medskip

Finally, we turn to Auslander-regular algebras.

Let $\Lambda$ be a $k$-algebra over a field $k$. Recall that
$\Lambda$ is called {\em Auslander-regular} if $\Lambda$ has finite
global dimension and satisfies the Gorenstein condition: if $p < q$
are non-negative integers and $M$ is a finitely generated (left or
right) $\Lambda$-module, then $\Ext^p_{\Lambda}(N, \Lambda) = 0$ for
every submodule N of $\Ext^q_{\Lambda\opp}(M, \Lambda)$. Here, if
$M$ is a right $\Lambda$-module, then $N$ is a left
$\Lambda$-module. Let $j(M)$ be the minimal number $r\ge 0$ such
that $\Ext^r_{\Lambda\opp}(M,\Lambda)\ne 0$. Then for any submodule
$N$ of $\Ext^{j(M)}_{\Lambda\opp}(M,\Lambda)$, we have
$\Ext^i_{\Lambda}(N,\Lambda)=0$ for $0<i<j(M).$ Thus:

\begin{Koro} Let $\Lambda$ be an Auslander-regular $k$-algebra, and
$M$ a finitely generated right $\Lambda$-module. Then, for any
submodule $X$ of $\Ext^{j(M)}_{\Lambda\opp}(M,\Lambda)$, and any
admissible subset $\Phi$ of ${\mathbb N}_{j(M)-2}$, the algebras
$\Ex{\Lambda}{\Phi}{\Lambda \oplus X}$ and
$\Ex{\Lambda}{\Phi}{\Lambda \oplus \Omega(X)}$ are derived
equivalent. \label{aregular} \end{Koro}

\subsection{Frobenius categories \label{frobenius}}

Let $\mathcal{A}$ be a {\em Frobenius} abelian category, that is,
${\cal A}$ is an abelian category with enough projective objects and
enough injective objects such that the projective objects coincides
with the injective objects. We  denote by $\underline{\cal A}$
the stable category of ${\cal A}$ modulo projective objects. It is
shown in \cite{HappelTriangle} that $\underline{\cal A}$ is a
triangulated category, in which the shift functor $[1]$ is just the
co-syzygy functor $\Omega^{-1}$, and the triangles in
$\underline{\cal A}$ are all induced by  short exact sequences in
${\cal A}$. For each morphism $f: U\ra V$ in ${\cal A}$, we denote
by $\underline{f}$ the image of $f$ under the canonical functor from
${\cal A}$ to $\underline{\cal A}$. Note that the objects of
$\underline{\cal A}$ are the same as those of $\cal A$.

\begin{Lem}
Let $\Phi$ be an admissible subset of $\mathbb{N}$, and let $M$,
$X$, and $Y$ be objects in ${\cal A}$. Then

$(1)$ For arbitrary $0\ne i\in \mathbb{N} $ and $U, U'\in {\cal A}$,
there is an isomorphism $$\Hom_{\Db{\cal A}}(U, U'[i])\simeq
\Hom_{\underline{\cal A}}(U, U'[i]), $$ which is functorial in $U$
and $U'$;

$(2)$ A monomorphism $\alpha: X\ra M_1$ in ${\cal A}$ is a left
$(\add(M),\Phi)$-approximation of $X$ in $\Db{\cal A}$ if and only
if $\underline{\alpha}$ is a left $(\add(M), \Phi)$-approximation of
$X$ in $\underline{\cal A}$;

$(3)$ An epimorphism $\beta: M_2\ra Y$ in ${\cal A}$ is a right
$(\add(M), -\Phi)$-approximation of $Y$ in $\Db{\cal A}$ if and only
if $\underline{\beta}$ is a right $(\add(M), -\Phi)$-approximation
of $Y$ in $\underline{\cal A}$. \label{lemFrobenius}
\end{Lem}

{\it Proof. } (1) For $0\ne i\in {\mathbb N}$, the isomorphisms
$$\Hom_{\Db{\cal A}}(U, U'[i])\simeq \Ext_{\cal A}^i(U, U')\simeq \Hom_{\underline{\cal A}}(U,
\Omega^{-i}U')=\Hom_{\underline{\cal A}}(U, U'[i]).$$ are functorial
in $U$ and $U'$. Thus (1) follows.

(2) First, let $0\ne i$ be in $\Phi$. By (1), there is a commutative
diagram
$$\xymatrix@C=15mm{
\Hom_{\Db{\cal A}}(M_1, M[i])\ar[r]^{(\alpha, M[i])}\ar[d]^{\simeq}
& \Hom_{\Db{\cal
A}}(X, M[i])\ar[d]^{\simeq}\\
\Hom_{\underline{\cal A}}(M_1, M[i])\ar[r]^{(\underline{\alpha},
M[i])} & \Hom_{\underline{\cal A}}(X, M[i]). }$$ Thus, the map
$\Hom_{\underline{\cal A}}(\underline{\alpha}, M[i])$ is surjective
if and only if $\Hom_{\Db{\cal A}}(\alpha, M[i])$ is surjective. Now
we consider the case $i=0$. If every morphism from $X$ to $M$ in
${\cal A}$ factorises through $\alpha$, then every morphism from $X$
to $M$ in $\underline{\cal A}$ factorises through
$\underline{\alpha}$. Conversely, assume that every morphism from
$X$ to $M$ in $\underline{\cal A}$  factorises through
$\underline{\alpha}$. Let $f: X\ra M$ be a morphism in ${\cal A}$.
Then $\underline{f}=\underline{\alpha}\underline{h}$ for some $h:
M_1\ra M$ in ${\cal A}$. Thus $f-\alpha h$ in ${\cal A}$ factorises
through a projective object $P$, say $f-\alpha h=st$ for some $s:
X\ra P$ and $t: P\ra M$ in ${\cal A}$. Since $P$ is also injective
and $\alpha$ is a monomorphism, there is some morphism $r: M_1\ra P$
such that $s=\alpha r$. Altogether, $f=\alpha h+st=\alpha h+\alpha
rt=\alpha(h+rt)$ factorises through $\alpha$. Thus the statement (2)
follows. The proof of (3) is similar to that of (2). $\square$

\begin{Prop}
Let $\Phi$ be an admissible subset of $\mathbb{N}$. Suppose that
${\cal A}$ is a Frobenius abelian category, that $M$ is an object in
${\cal A}$, and that $0\ra X\raf{\alpha} M_1\raf{\beta} Y\ra 0$ is a
short exact sequence in ${\cal A}$ with $M_1\in\add(M)$ such that
the induced triangle $X\lraf{\underline{\alpha}}
M_1\lraf{\underline{\beta}} Y\lra X[1]$ in $\underline{\cal A}$
satisfies the conditions in {\rm Theorem \ref{verygthm1}}.  Then the
algebras $\Ex{\D{\cal A}}{\Phi}{M\oplus Y}$ and $\Ex{\D{\cal
A}}{\Phi}{X\oplus M}$ are derived equivalent. \label{FrobProp}
\end{Prop}

{\it Proof.} This follows from Lemma \ref{lemFrobenius} and
Proposition \ref{vanish2}. $\square$

\begin{Koro}
Suppose that ${\cal A}$ is a Frobenius abelian category and $M$ is
an object in ${\cal A}$. Let $0\ra X\ra M_1\ra Y\ra 0$ with
$M_1\in\add(M)$ be a short exact sequence in ${\cal A}$ such that
the induced triangle in $\underline{\cal A}$ is an $\add(M)$-split
triangle.  Then $\End_{\cal A}(M\oplus Y)$ and $\End_{\cal
A}(X\oplus M)$ are derived equivalent. \label{FrobCor1}
\end{Koro}

{\it Proof.} Taking $\Phi:=\{0\}$, the corollary follows from
Proposition \ref{FrobProp}. $\square$

\medskip
{\it Remark}. If $\cal A$ is a Frobenius (not necessarily abelian)
category, then Corollary \ref{FrobCor1} is still true. For the
precise definition of a Fronenius category, we refer the reader to
\cite{HappelTriangle}.

\medskip
The module category of a self-injective Artin algebra is a Frobenius
abelian category. In this case, we have the following corollary.

\begin{Koro}
Let $A$ be a self-injective algebra, and let $M$ be an $A$-module.
Suppose $X\ra M_1\ra Y\ra X[1]$ is an  $\add(M)$-split triangle in
$\stmodcat{A}$. Then $\End_A(A\oplus M\oplus X)$ and $\End_A(A\oplus
M\oplus Y)$ are derived equivalent.
\end{Koro}

{\it Proof.} Since all triangles in $\stmodcat{A}$ are induced by
short exact sequences in $\modcat{A}$, there is a short exact
sequence $0\ra X\ra M_1\oplus P\ra Y\ra 0$ in $A$-mod with $P$
projective such that the induced triangle is isomorphic to the given
triangle $X\ra M_1\ra Y\ra X[1]$ in $\stmodcat{A}$. The triangle
$X\ra M_1\oplus P\ra Y\ra X[1]$ also is an  $\add(A\oplus M)$-split
triangle in $\stmodcat{A}$. The corollary then follows from
Corollary \ref{FrobCor1}. $\square$

\subsection{Calabi-Yau categories}

The theory of Calabi-Yau and cluster categories provides very
natural contexts for our construction of derived equivalences.

Let $k$ be a field, and let ${\cal T}$ be a $k$-linear triangulated
category which is Hom-finite, that is, the Hom-space $\Hom_{\cal
T}(X, Y)$ is finite dimensional over $k$ for all $X$ and $Y$ in
${\cal T}$.

Recall that ${\cal T}$ is called \emph{$(n+1)$-Calabi-Yau} for some
non-negative integer $n$ if there is a natural isomorphism between
$D\Hom_{\cal T}(X,Y)$ and $\Hom_{\cal T}(Y, X[n+1])$ for all $X$ and
$Y$ in ${\cal T}$, where $D=\Hom_k(-, k)$ is the usual duality. It
follows that ${\mathscr X}_{\cal T}^n(M)= {\mathscr Y}_{\cal
T}^n(M)$ for $M\in {\cal T}$. (See \cite{kellercy} for more
information on Calabai-Yau categories.)

Note that if $\Phi = \{0, 1, \cdots, n\}$, then $n-i\in\Phi$ for
each $i\in\Phi$.

\begin{Lem}
Let $\Phi=\{0, 1,\cdots, n\}$. Suppose that ${\cal T}$ is an
$(n+1)$-Calabi-Yau triangulated category, and that $M$ is an object
in ${\cal T}$. Let $X\lraf{\alpha} M_1\lraf{\beta} Y\lra X[1]$ be a
triangle in ${\cal T}$ with $M_1\in\add(M)$. Then:

$(1)$ The morphism $\alpha$ is a left $(\add(M),
\Phi)$-approximation of $X$ if and only if the morphism $\beta$ is a
right $(\add(M), -\Phi)$-approximation of $Y$;

$(2)$ If $\alpha$ is a left $(\add(M), \Phi)$-approximation of $X$
and if $M$ is $n$-self-orthogonal, then $X\in{\mathscr
X}^{n}(M)\cap{\mathscr Y}^{n}(M)$ and $Y\in{\mathscr
X}^n(M)\cap{\mathscr Y}^n(M)$. \label{propCY}
\end{Lem}

{\it Proof.} We will abbreviate $\Hom_{\cal T}(-, -)$ by $(-, -)$.
First we assume that $\alpha$ is a left $(\add(M),
\Phi)$-approximation of $X$. Now, for each $i\in\Phi$, there is a
commutative diagram with exact rows
$$\xymatrix@C=17mm{
(M[-i], M_1)\ar[r]^{(M[-i], \beta)}\ar[d]^{\simeq} & (M[-i], Y)\ar[d]^{\simeq}\\
(M, M_1[i])\ar[r]^{(M, \beta[i])} &  (M, Y[i])\ar[r] &
(M, X[i+1])\ar[r]^{(M, \alpha[i+1])}\ar[d]^{\simeq} & (M, M_1[i+1])\ar[d]^{\simeq}\\
&& D(X, M[n-i])\ar[r]^{D(\alpha,M[n-i])} & D(M_1, M[n-i]). }$$ Since
$n-i$ is in $\Phi$, and since $\alpha$ is a left $(\add(M),
\Phi)$-approximation of $X$, the map $(\alpha, M[n-i])$ is
surjective, and consequently $D(\alpha, M[n-i])$ is injective. Hence
$(M, \alpha[i+1])$ is injective, and therefore $(M[-i], \beta)$ is
surjective. This shows that $\beta$ is a right $(\add(M),
-\Phi)$-approximation of $Y$. The proof of the other implication in
(1) can be done similarly.

(2) It follows from (1) and the comment before Corollary
\ref{cortoverygthm1} that $X\in{\mathscr X}_{\cal T}^{\Phi}(M)$ and
$Y\in{\mathscr Y}_{\cal T}^{\Phi}(M)$. Since ${\cal T}$ is
$(n+1)$-Calabi-Yau, we have $(M, X[i])\simeq D(X, M[n+1-i])=0$, and
$(M, Y[i])\simeq D(Y, M[n+1-i])=0$ for all $0\neq i\in\Phi$. Thus
$X\in{\mathscr Y}_{\cal T}^{\Phi}(M)$ and $Y\in{\mathscr X}_{\cal
T}^{\Phi}(M)$. $\square$

\begin{Koro}
Let $\Phi=\{0, 1, \cdots, n\}$, and let ${\cal T}$ be an
$(n+1)$-Calabi-Yau triangulated category. Suppose that $M$ is
$n$-self-orthogonal and $Y\in{\mathscr Y}^n(M)$. Let $X\lraf{\alpha}
M_1\lraf{\beta} Y\lraf{w} X[1]$ be a triangle in ${\cal T}$ with
$\beta$ a right $\add(M)$-approximation of $Y$.  Then the algebras
$\Ex{\cal T}{\Phi}{M\oplus X}/I$ and $\Ex{\cal T}{\Phi}{M\oplus
Y}/J$ are derived equivalent, where $I$ and $J$ are defined as in
{\rm Theorem \ref{verygthm1}}. \label{cycor1}
\end{Koro}

{\it Proof.} Since $Y\in{\mathscr Y}_{\cal T}^{\Phi}(M)$, for each
$0\neq i\in\Phi$, the map $(M[-i], M_1)\lra (M[-i], Y)=0$ induced by
$\beta$ is surjective. Taking into account that $\beta$ is a right
$\add(M)$-approximation of $Y$, we see that $\beta$ is, in fact,  a
right $(\add(M), -\Phi)$-approximation of $Y$. By Proposition
\ref{propCY} (1), the map $\alpha$ is a left $(\add(M),
\Phi)$-approximation of $X$. Since $M$ is $n$-self-orthogonal, the
proof can be finished by applying Proposition \ref{propCY} (2) and
Corollary \ref{cortoverygthm1} to the triangle. $\square$

\medskip
Corollary \ref{cycor1} is related to mutations in a Calabi-Yau
category. Here are some definitions from \cite{IY}.

Let $\cal T$ be an $(n+1)$-Calabi-Yau category. An object $T$ in
$\cal T$ is called an \emph{$n$-cluster tilting object} if $T$ is
$n$-self-orthogonal, and if any $X \in {\cal T}$ with $\Ext^i_{\cal
T}(T,X)=0$ for $1\le i\le n$ is in $\add(T)$. The object $T$ is
called {\em basic} if the multiplicity of each indecomposable direct
summand of $T$ is one.

Let $T$ be an  $n$-cluster basic tilting object in an
$(n+1)$-Calabi-Yau category $\cal T$, and $Y$ a direct summand of
$T$, that is, $T = Y\oplus M$. Let $\beta: M_1\ra Y$ be a minimal
right $\add(M)$-approximation of $Y$, and let
$$X\lraf{\alpha}M_1\lraf{\beta}Y\lra X[1]
$$ be a triangle containing $\beta$. Note that we allow $Y$ to be decomposable, and
that $X$ is indecomposable if and only if $Y$ is indecomposable. The
object $X\oplus M$ is called the \emph{left mutation} of $T$ at $Y$.
In the case of tilting modules, $X$ is called a \emph{tilting
complement} to $M$ in the literature (see, for example, \cite{HU}).
It was pointed out in \cite{IY} that the left mutation of $T$ at $Y$
is again an $n$-cluster tilting object (for some special cases, see
\cite{Bondal, Gor}, and also \cite[p.314]{LM}). In fact, this can be
seen in the following way: The proof of Corollary \ref{cycor1} and
comment (b) on the conditions of Theorem \ref{verygthm1} imply that
$T':= M\oplus X$ is $n$-self-orthogonal. Morover, let $X'\in
{\mathscr X}^n(T')$ and consider a triangle $X'\lraf{\alpha'} M'\ra
Y'\ra X'[1]$ with $\alpha'$ a left $\add(M)$-approximation of $X'$.
Then $Y'\in {\mathscr X}^n(T)$ by Lemma \ref{propCY} and the comment
(b). Thus $Y'\in \add(T)$, $X'\in \add(T')$, and $T':=X\oplus M$ is
again an $n$-cluster tilting object in $\cal T$. The notion of a
right mutation of $T$ at $Y$ is dual.

Usually, $\End_{\cal T}(X\oplus M)$ and $\End_{\cal T}(M\oplus Y)$
are not derived equivalent. When they are derived equivalent may be
an interesting question. Here is a sufficient condition.

\begin{Koro} Let $\Lambda :=\End_{\cal T}(X\oplus M)$ and $\Gamma :=\End_{\cal T}(M\oplus
Y)$. Then

$(1)$ $\End_{\cal T}(X\oplus M)/I$ and $\End_{\cal T}(M\oplus Y)/J$
are derived equivalent.

$(2)$ Suppose that $Y$ is indecomposable. Let $S_X$ be the simple
$\Lambda$-module corresponding to $X$, and let $S_Y$ be the simple
$\Gamma$-module corresponding to $Y$. Suppose that $S_Y$ is not a
submodule of $\Gamma$, and $S_X$ is not a quotient of $D(\Lambda)$.
Then $\Lambda$ and $\Gamma$ are derived equivalent. \label{mutation}
\end{Koro}

{\it Proof.} Statement (1) is a direct consequence of Corollary
\ref{cycor1}, and (2) follows from (1) and Proposition
\ref{bigideal}. $\square$

\medskip
{\it Remark.} Consider a $2$-Calabi-Yau category, and assume that
$\Ext^1_{\Gamma}(S_Y,S_Y)=0$. Then we re-obtain the result
\cite[Theorem 5.3]{ladkani} from Corollary \ref{mutation} (2).

\section{Examples \label{sect5}}
First, we present an explicit example which satisfies all conditions
in Theorem \ref{verygthm1}.

\medskip
{\parindent=0pt\bf Example 1.} Let $k$ be an algebraically closed
field of characteristic $2$, and let $A:=kA_4$ be the group algebra
of the alternating group $A_4$. Then there are three simple
$A$-modules, which are denoted $k, \omega$, and $\bar{\omega}$,
respectively. Their projective covers are $P(k)$, $P(\omega)$ and
$P(\bar{\omega})$, respectively. It is well known that $kA_4$ is
Morita equivalent to the following algebra given by quiver

$$
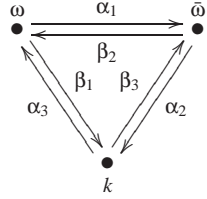
\xymatrix@R=6mm@C=8mm{
 \bullet\ar@<2pt>[rr]^{\alpha_1}^(0){\omega}^(1){\bar{\omega}}\ar@<2pt>[rdd]^{\beta_1}_(1.1){k}
 & & \bullet\ar@<2pt>[ldd]^{\alpha_2}\ar@<2pt>[ll]^{\beta_2}\\ & & \\
 &\bullet\ar@<2pt>[luu]^{\alpha_3}\ar@<2pt>[ruu]^{\beta_3} &
} $$ and relations
$\alpha_i\beta_{i+1}-\beta_i\alpha_{i+2}=\alpha_i\alpha_{i+1}=\beta_i\beta_{i-1}=0$,
where the subscripts are considered modulo $3$.

As this algebra is symmetric, the Auslander-Reiten translation
$D\Tr$ is just the second syzygy $\Omega^2$. The Auslander-Reiten
quiver of this algebra is well-known to have a component of the
following form:

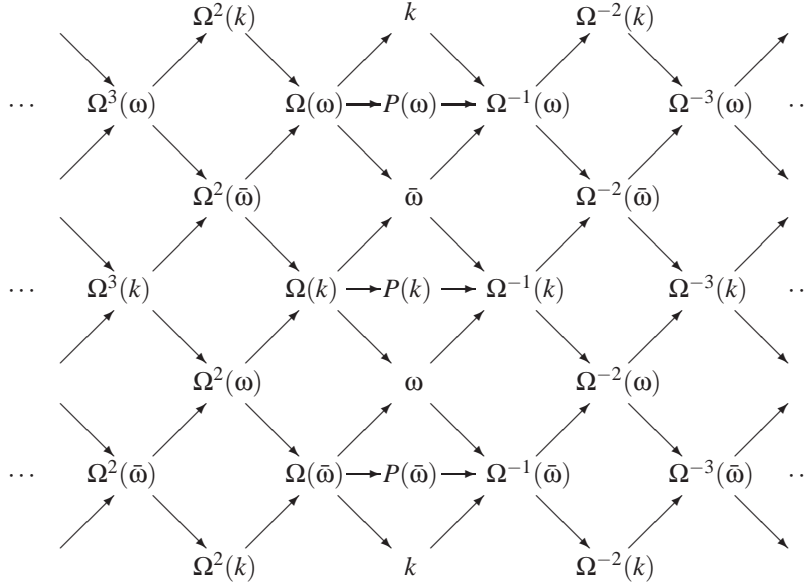
\begin{center}
{\small
\begin{picture}(330, 260)
    % \put(40, 315){\includegraphics{longdots}}
     %\put(322, 315){\includegraphics{longdots}}
     %\put(135, 315){\includegraphics{longdots}}
     %\put(179, 315){\includegraphics{longdots}}

    \put(105, 230){\vector(1,-1){20}}
    \put(70, 210){\vector(1,1){20}}
    \put(85, 233){$\Omega ^{2} (k)$}
    \put(70, 195){\vector(1,-1){20}}
    \put(105, 175){\vector(1,1){20}}
    \put(85, 165){$\Omega ^{2} (\bar{\omega})$}
    \put(70, 140){\vector(1,1){20}}
    \put(105, 160){\vector(1,-1){20}}
    \put(70, 125){\vector(1,-1){20}}
    \put(105, 105){\vector(1,1){20}}
    \put(85, 95){$\Omega ^{2} (\omega)$}
    \put(70, 70){\vector(1,1){20}}
    \put(105, 90){\vector(1,-1){20}}
    \put(70, 55){\vector(1,-1){20}}
    \put(105, 35){\vector(1,1){20}}
    \put(85, 25){$\Omega ^{2} (k)$}
   % \put(135, 47){\includegraphics{longdots}}
    \put(15, 200){$\cdots$}
    \put(45, 200){$\Omega ^{3} (\omega)$}
    \put(15, 130){$\cdots$}
    \put(45, 130){$\Omega ^{3} (k)$}
    \put(15, 60){$\cdots$}
    \put(45, 60){$\Omega ^{2} (\bar{\omega})$}

    \put(215, 210){\vector(1,1){20}}
    \put(250, 230){\vector(1,-1){20}}
    \put(230, 233){$\Omega ^{-2} (k)$}
    \put(290, 210){\vector(1,1){20}}
    \put(290, 140){\vector(1,1){20}}
    \put(290, 70){\vector(1,1){20}}
    \put(290, 195){\vector(1,-1){20}}
    \put(290, 125){\vector(1,-1){20}}
    \put(290, 55){\vector(1,-1){20}}

    \put(35, 230){\vector(1,-1){20}}
    \put(35, 160){\vector(1,-1){20}}
    \put(35, 90){\vector(1,-1){20}}
    \put(35, 175){\vector(1,1){20}}
    \put(35, 105){\vector(1,1){20}}
    \put(35, 35){\vector(1,1){20}}

    \put(215, 195){\vector(1,-1){20}}
    \put(250, 175){\vector(1,1){20}}
    \put(230, 165){$\Omega ^{-2} (\bar{\omega})$}
    \put(215, 140){\vector(1,1){20}}
    \put(250, 160){\vector(1,-1){20}}
    \put(215, 125){\vector(1,-1){20}}
    \put(250, 105){\vector(1,1){20}}
    \put(230, 95){$\Omega ^{-2} (\omega)$}
    \put(215, 70){\vector(1,1){20}}
    \put(250, 90){\vector(1,-1){20}}
    \put(215, 55){\vector(1,-1){20}}
    \put(250, 35){\vector(1,1){20}}
    \put(230, 25){$\Omega ^{-2} (k)$}
    %\put(322, 47){\includegraphics{longdots}}
    \put(310, 200){$\cdots$}
    \put(265, 200){$\Omega ^{-3} (\omega)$}
    \put(310, 130){$\cdots$}
    \put(265, 130){$\Omega ^{-3} (k)$}
    \put(310, 60){$\cdots$}
    \put(265, 60){$\Omega ^{-3} (\bar{\omega})$}

    \put(165, 235){$k$}
    \put(140, 210){\vector(1,1){20}}
    \put(175, 230){\vector(1,-1){20}}

   \put(120, 200){$\Omega (\omega)$}
   \put(143, 203){\vector(1,0){13}}
   \put(157, 200){$P(\omega)$}
   \put(179, 203){\vector(1,0){13}}
   \put(196, 200){$\Omega ^{-1} (\omega)$}

   \put(140, 195){\vector(1,-1){20}}
   \put(175, 175){\vector(1,1){20}}

   \put(165, 165){$\bar{\omega}$}
   \put(140, 140){\vector(1,1){20}}
   \put(175, 160){\vector(1,-1){20}}

   \put(120, 130){$\Omega (k)$}
   \put(143, 133){\vector(1,0){13}}
   \put(157, 130){$P(k)$}
   \put(179, 133){\vector(1,0){13}}
   \put(196, 130){$\Omega ^{-1} (k)$}

   \put(140, 125){\vector(1,-1){20}}
   \put(175, 105){\vector(1,1){20}}

   \put(165, 95){$\omega$}
   \put(140, 70){\vector(1,1){20}}
   \put(175, 90){\vector(1,-1){20}}

   \put(120, 60){$\Omega (\bar{\omega})$}
   \put(143, 63){\vector(1,0){13}}
   \put(157, 60){$P(\bar{\omega})$}
   \put(179, 63){\vector(1,0){13}}
   \put(196, 60){$\Omega ^{-1} (\bar{\omega})$}

   \put(140, 55){\vector(1,-1){20}}
   \put(175, 35){\vector(1,1){20}}
   \put(165, 25){$k$}
 % \put(179, 47){\includegraphics{longdots}}
 % \put(40, 47){\includegraphics{longdots}}

%\put(115, 20){\line(1, 0){110}}\put(115, 20){\line(0, 1){225}}\put(115, 245){\line(1, 0){110}}\put(225, 20){\line(0, 1){225}}

%\iffalse \put(79, 10){\line(1,0){180}}\put(79,10){\line(0,1){245}}\put(79, 255){\line(1,0){180}}\put(259,10){\line(0,1){245}}\fi
\end{picture}}
\end{center}

Consider the Auslander-Reiten sequence $$0\lra
\Omega^{3}(\omega)\lra\Omega^{2}(k)\oplus
\Omega^{2}(\bar{\omega})\lra \Omega(\omega)\lra 0.$$ Let
$X=\Omega^{3}(\omega)$, $Y=\Omega(\omega)$, and
$M=\Omega^{2}(k)\oplus \Omega^{2}(\bar{\omega})$. This sequence
provides an Auslander-Reiten triangle in the triangulated category
$\stmodcat{A}$:

$$X\lra M\lra Y\lra X[1].$$
We shall check that this triangle satisfies the conditions of
Theorem \ref{verygthm1}.

We choose $\Phi=\{0,1\}$ and $F=[1]$. Since this is an
Auslander-Reiten triangle in $\stmodcat{A}$, the map $X\ra M$ is a
left $(\add(M), \Phi)$-approximation of $X$, and the map $M\ra Y$ is a right
$(\add(M), -\Phi)$-approximation of $Y$ (see the example at the end
of Section \ref{sect2}). It follows from the above Auslander-Reiten
quiver of $A$ that $\Ext^1_A(M,X) \simeq
\StHom_A(M,\Omega^{-1}(X))\simeq \StHom_A(\Omega^{2}(k)\oplus
\Omega^{2}(\bar{\omega}), \Omega^{2}(\omega))=0$ and
$\Ext^1_A(Y,M)\simeq
\StHom_A(Y,\Omega^{-1}(M))=\StHom_A(\Omega(\omega), \Omega(k)\oplus
\Omega(\bar{\omega}))=0.$ Thus the above triangle in $\stmodcat{A}$
satisfies all conditions in Theorem \ref{verygthm1}, and therefore,
by Proposition \ref{FrobProp}, the algebras $\Ex{A}{\Phi}{M\oplus
X}$ and $\Ex{A}{\Phi}{M\oplus Y}$ are derived equivalent.

Furthermore, we have $\Ext_A^1(M, M)\simeq \StHom_A(M,
\Omega^{-1}M)\simeq \StHom_A(\Omega(k)\oplus \Omega(\bar{\omega}),
k\oplus \bar{\omega})$. There is an epimorphism from $\Omega(k)$ to
$\bar{\omega}$ and an epimorphism from $\Omega(\bar{\omega})$ to
$k$. The latter cannot factorise through a projective module, we get
dim$_k\Ext_A^1(M, M)=2$. Moreover, there is an epimorphism from
$\Omega(k)$ to $\omega$ and an epimorphism from
$\Omega(\bar{\omega})$ to $\omega$. This implies dim$_k\Ext_A^1(M,
Y)=2$. Similarly, dim$_k\Ext_A^1(X, M)=2$. Note that all the
indecomposable modules appearing in the Auslander-Reiten triangle
are $1$-self-orthogonal. A more precise calculation shows that
$\dim_k\Ex{A}{\Phi}{M\oplus X}=33$ and $\dim_k\Ex{A}{\Phi}{M\oplus
Y}=21.$

\medskip
The following example shows that the Ext-orthogonality conditions in
Corollary \ref{thm1adss} and therefore in Theorem \ref{verygthm1}
cannot be dropped.

\medskip
{\parindent=0pt\bf Example 2.} Let $A$ be the algebra (over a field
$k$) given by the following quiver with relations:

{\unitlength=0.5cm \special{em:linewidth 0.4pt}
\linethickness{0.4pt}
\begin{center}
\begin{picture}(20,3.5)
\put(5,1.3){$1$}\put(7.3,1.3){$2$}
\put(5,2){$\bullet$}\put(7.3,2.2){\vector(-1,0){1.8}}\put(7.5,2){$\bullet$}\put(8.7,2.3){\circle{2}}
\put(9.7,2.2){\vector(-1,-2){0.2}}

\put(6.5, 2.5){$\beta$}\put(10, 2){$\alpha ,$}\put(11.5,2){$\alpha^2
= 0 = \alpha\beta.$}
\end{picture}
\end{center} }\vspace{-0.5cm}
This example is in a class of examples constructed by Small
\cite{Small}. The algebra $A$ is of finite representation type, its
finitistic dimension equals one, while the finitistic dimension of
the opposite algebra $A\opp$ is zero.

We denote by $S(i)$ and $P(i)$ the simple and projective modules
corresponding to the vertex $i$, respectively. Let $M_i$ be the
quotient module of $P(2)$ by $S(i)$, and $M:=M_1\oplus M_2=D(A_A)$,
where $D$ is the usual duality. Then there is an Auslander-Reiten
sequence
$$0\lra X:=P(2)\lra M\lra S(2)=: Y\lra 0. $$
This is an $\add(M)$-split sequence in $A$-mod.

If we take $\Phi=\{0,1\}$, then $\Ex{A}{\Phi}{X\oplus
M}=\End_A(X\oplus M)$. An easy calculation shows that $
\End_A(X\oplus M)$ is a quasi-hereditary algebra, and thus has
finite global dimension. The algebra $\Ex{A}{\Phi}{M\oplus Y}$
contains a loop which is given by the short exact sequence induced
by the loop $\alpha$ at the vertex $2$. Thus it has infinite global
dimension by the 'no loops theorem'. It follows that
$\Ex{A}{\Phi}{X\oplus M}$ and $\Ex{A}{\Phi}{M\oplus Y}$ cannot be
derived equivalent since derived equivalences preserve the
finiteness of global dimensions. Also, one can see that
$\Ext_A^i(X,M)=0=\Ext^1_A(M,X)$ and $\Ext^i_A(Y,M)=0\ne
\Ext^1_A(M,Y)$ for $i\ge 1$. This example shows that the
orthogonality conditions in Corollary \ref{thm1adss} cannot be
omitted. Moreover, it shows that the result in \cite[Theorem
1.1]{hx2} cannot be extended from endomorphism algebras to
$\Phi$-Yoneda algebras without any additional conditions.

\bigskip
\begin{appendix}
\section{ A two functors version of Theorem
\ref{thm1}}

\bigskip
\noindent In Theorem \ref{verygthm1}, there is only one functor $F$
involved. When working with the derived category of a
hereditary algebra, or the stable category of a self-injective
algebra, or the derived category of coherent sheaves of a projective
variety over $\mathbb C$, 
%(see \cite{LM}, \cite{st})
apart from the
shift functor there are other prominent functors, for example the
Auslander-Reiten translation DTr. To have available a general statement of
construction of derived equivalences, which is similar to Theorem
\ref{verygthm1}, we define $\Phi$-perforated Yoneda algebras
for two functors over a triangulated category, and formulate a
two-functor version of Theorem \ref{verygthm1}. In this appendix, we
summarise the ingredients for a possible generalisation of Theorem
\ref{verygthm1}. The proof of this generalisation is
analogous to that of Theorem \ref{verygthm1}, but more technical
and tedious. So we omit it here.

Let $\Phi$ be a subset of ${\mathbb N} \times {\mathbb N}$, which we
consider as a semigroup with ordinary addition. Let $\cal T$ be a
triangulated $R$-category with shift functor [1], and let $X$ be
an object in $\cal T$.

Suppose that $F$ and $G$ are two triangle functors from $\cal T$ to
itself, such that $FG$ is naturally isomorphic to $GF$. For $X$ in
$\cal T$, let $\delta(i,j,X): F^jG^iX\ra G^iF^jX$ be an isomorphism
induced from the natural transformation $FG\sim GF$. Then we define

$$\Ex{\cal T}{F,G,\Phi}{X} := \displaystyle \bigoplus_{(i,j)\in
\Phi} \Hom_{\cal T}(X, G^iF^jX),$$ with elements of the form
$(f_{i,j})_{(i,j)\in \Phi}$, where $f_{i,j}: X\ra
G^iF^jX$. The multiplication on $\Ex{\cal T}{F,G,\Phi}{X}$ is given by 

$$ \Big(f_{i,j}\Big)_{(i,j)\in \Phi}\cdot \Big(g_{i,j}\Big)_{(i,j)\in \Phi}
= \Big(\sum_{{{(p,q),(u, v)\in\Phi}\atop {(u+p,
v+q)=(l,m)\in\Phi}}}f_{u,v}(G^uF^vg_{p,q})(G^u\delta(p,v,F^qX)\Big)_{(l,m)\in
\Phi\times\Psi}.$$ 
If $F$ and $G$ are invertible, then $\Phi$ can be chosen a
subset of ${\mathbb Z} \times {\mathbb Z}$.

A general model for the above definition is: Given a bi-graded
algebra $\Lambda = \bigoplus_{i,j\in {\mathbb Z}} \Lambda_{i,j}$, we
define $\Lambda(\Phi)= \displaystyle\bigoplus_{(i,j)\in \Phi}
\Lambda_{i,j}$, and a multiplication by $ a_{i,j}\cdot a_{p,q}=
a_{i,j}a_{p,q}$ if $(i+p, j+q) \in \Phi$, and zero otherwise. If
$\Phi$ is admissible, for example, $\Phi$ is the cartesian product of
two admissible sets in $\mathbb Z$, then $\Lambda(\Phi)$ is an
associative algebra. So, we have to check that $\Ex{\cal T}{F,G}{X}
:= \displaystyle \bigoplus_{i,j\in {\mathbb Z}} \Hom_{\cal T}(X,
G^iF^jX)$ is an associative algebra with respect to the above
multiplication. This can be based on the
following lemma.

\begin{Lem} Suppose $F$ and $G$ are two triangle functors from $\cal T$ to
itself such that $FG$ is naturally isomorphic to $GF$. For any
triangle functor $L$ from $\cal T$ to itself, there is a natural
isomorphism $\delta(i,j, L): F^jG^iL\lra G^iF^jL$ for all $i,j\ge 0$
such that, for $p, q, r, s\in\mathbb{N}$,

\smallskip
$(1)$ $\delta(p+q, r, L)=\delta(p, r, G^{q}L) \big(G^p\delta(q, r,
L)\big);$

\medskip
$(2)$ $\delta(p, r+s, L)=\big(F^{s}\delta(p, r, L)\big)\delta(p, s,
F^rL)$. \label{LemAsso}

\end{Lem}

{\it Proof.} For functors $L_1$ and $L_2$ from ${\cal T}$ to itself,
we define $L_1\delta(1, 1, L_2): L_1FGL_2\ra L_1GFL_2$ to be the
induced natural isomorphism from the functor $L_1FGL_2$ to the
functor $L_1GFL_2$. So, $\delta(1,1,1_{\cal T})$ is just the
given natural isomorphism from $FG$ to $GF$. Now we shall construct
inductively a natural isomorphism $\delta(i,j, L)$ from $F^jG^iL$ to
$G^iF^jL$ for all non-negative integers $i$ and $j$ and functors $L$
from $\cal T$ to itself.

If $i=0$ or $j=0$, then $F^jG^iL=G^iF^jL$, and we define
$\delta(i,j, L)$ to be the identity natural transformation. For each
positive integer $j>1$, we assume that $\delta(1, j-1, L)$ is
defined. Now we define
$$\delta(1, j, L):=\big(F\delta(1,j-1,L)\big)\delta(1,1,F^{j-1}L).$$
For each positive integer $i>1$, assume that $\delta(i-1, j, L)$ is
defined. We define
$$\delta(i,j, L):=\delta(1,j, G^{i-1}L)\big(G\delta(i-1,j, L)\big).$$

(1) It is straightforward to check that (1) holds for $p+q\leq 2$.
We shall prove (1) by induction on $p+q$. Now assume that $p+q>2$.
Then we have
$$\begin{array}{rl}
 \delta(p+q, r, L) & = \delta(1,r, G^{p+q-1}L)\big(G\delta(p+q-1,r, L)\big)\quad \mbox{(by definition)}\\
 & =  \delta(1,r, G^{p+q-1}L) G\Big(\delta(p-1, r, G^{q}L) \big(G^{p-1}\delta(q, r,
    L)\big)\Big) \quad \mbox{(by induction)}\\
    &= \Big(\delta(1,r, G^{p+q-1}L)\big(G\delta(p-1, r,
    G^{q}L)\big)\Big)\big(G^{p}\delta(q, r,
    L)\big)\\
    & =\delta(p, r, G^qL)\big(G^p\delta(q, r, L)\big) \quad \mbox{(by definition)}.
\end{array}$$
This proves (1).

(2) We first prove (2) for $p=0, 1$. If $p=0$, then (2) is clearly
true. Now suppose $p=1$. We shall show (2) by induction on $r+s$. In
fact, if $r+s\leq 2$, it is straightforward to check (2). Now we
assume that $r+s>2$. Then we have
$$\begin{array}{rl}
 \delta(1, r+s, L) & = \big(F\delta(1,r+s-1,L)\big)\delta(1,1,F^{r+s-1}L)\quad \mbox{(by definition)}\\
 & = F\Big(\big(F^{s-1}\delta(1, r, L)\big)\delta(1, s-1,
    F^rL)\Big)\delta(1,1,F^{r+s-1}L) \quad \mbox{(by induction)}\\
 & = \big(F^s\delta(1, r, L)\big) \Big(\big(F\delta(1, s-1,
    F^rL)\big)\delta(1,1,F^{r+s-1}L)\Big)\\
 & = \big(F^s\delta(1, r, L)\big)\delta(1, s, F^rL)\quad \mbox{(by definition)}.
 \end{array}$$
This proves (2) for $p=1$. Now assume $p>1$. Then {\small
$$\begin{array}{rcl}
 \delta(p, r+s, L) &= & \delta(1, r+s, G^{p-1}L) \big(G\delta(p-1, r+s, L)\big)\quad \mbox{(by definition)}\\
 & = & \big(F^s\delta(1, r, G^{p-1}L)\big)\delta(1, s, F^rG^{p-1}L)G\Big(\big(F^s\delta(p-1, r, L)\big)\delta(p-1, s,
F^rL)\Big)\quad \mbox{(by induction)} \\
& =& \big(F^s\delta(1, r, G^{p-1}L)\big) \Big(\delta(1, s,
F^rG^{p-1}L)\big(GF^s\delta(p-1, r, L)\big)\Big)\big(G\delta(p-1, s,
F^rL)\big).
\end{array}$$}
Since $\delta(1, s, F^rG^{p-1}L)$ is a natural transformation from
$F^sGF^rG^{p-1}L$ to $GF^s F^rG^{p-1}L$, the following
diagram of natural transformations is commutative:
$$\xymatrix@C=30mm{
  F^sGF^rG^{p-1}L \ar[r]^{\delta(1, s,
F^rG^{p-1}L)}\ar[d]_{F^sG\delta(p-1, r, L)} &
GF^sF^rG^{p-1}L\ar[d]^{GF^s\delta(p-1, r, L)}\\
  F^sGG^{p-1}F^rL \ar[r]^{\delta(1, s,
G^{p-1}F^rL)} & GF^sG^{p-1}F^rL.\\
}$$ Hence
$$\begin{array}{rcl}
 \delta(p, r+s, L) &= & \big(F^s\delta(1, r, G^{p-1}L)\big) \Big(\delta(1, s,
F^rG^{p-1}L)\big(GF^s\delta(p-1, r, L)\big)\Big)\big(G\delta(p-1, s, F^rL)\big)\\
&=& \big(F^s\delta(1, r, G^{p-1}L)\big) \Big(\big(F^sG\delta(p-1, r,
L)\big)\delta(1, s, G^{p-1}F^rL)\Big)\big(G\delta(p-1, s,
F^rL)\big)\\
& =& F^s\Big(\delta(1, r, G^{p-1}L)\big(G\delta(p-1, r,
L)\big)\Big)\Big(\delta(1, s, G^{p-1}F^rL)\big(G\delta(p-1, s,
F^rL)\big)\Big)\\
& =& \big(F^{s}\delta(p, r, L)\big)\delta(p, s,
    F^rL).
\end{array}$$
This proves (2). $\square$

\medskip
{\it Remark.} If, in addition, $F$ and $G$ are invertible, then
Lemma \ref{LemAsso} remains valid for $i,j, p,q,r$ and $s$ any integers.

\medskip
Let $\cal D$ be a full subcategory of $\cal T$, and $X$ an object of
$\cal T$. A morphism $f: X\ra D$ with $D\in \cal D$ is called a left
$({\cal D},F,G,\Phi)$-approximation of $X$ if $\Hom_{\cal
T}(f,G^iF^jD'): \Hom_{\cal T}(D,G^iF^jD')\lra \Hom_{\cal
T}(X,G^iF^jD')$ is surjective for every object $D'\in {\cal D}$ and
$(i,j)\in \Phi$. Dually, we define the right $({\cal
D},F,G,\Phi)$-approximation of $X$.

Given a triangle $0\ra X\stackrel{\alpha}{\lra}
M_1\stackrel{\beta}{\lra} Y \stackrel{w}{\lra} X[1]$ in $\mathcal T$
with $M_1\in \add(M)$ for a fixed $M\in {\cal T}$, we define
$\widetilde{w}[-1]= (-w[-1],0): Y[-1]\ra X\oplus M$, $\bar{w}=(0,
w)^T$, where $(0, w)^T$ stands for the transpose of the matrix
$(0,w)$, and

$$\begin{array}{rl} I  &:= \{x=(x_{i,j})\in \Ex{\cal T}{F, G, \Phi}{X\oplus M}\mid x_{i,j}=0
\; \mbox{for\;} (0,0)\ne (i,j)\in \Phi, \;\mbox{and} \\
&  \hspace{0.7cm} x_{0,0} \;\mbox{ factors through}\;
\add(M)\;\mbox{and}\; \widetilde{w}[-1] \},\end{array}$$
$$\begin{array}{rl}J  &:= \{ y=(y_{i,j})\in \Ex{\cal T}{F,G,\Phi}{M\oplus Y}\mid y_{i,j}=0
\;\mbox{for\;} (0,0)\ne (i,j)\in \Phi, \;\mbox{and}\\ &
\hspace{0.7cm} y_{0,0} \;\mbox{factors through}\; \add(M)\;
\mbox{and}\; \bar{w} \}.\end{array}$$

Now, Theorem \ref{thm1} generalises as follows:

\begin{Theo} Let $\Phi$ be an admissible subset of ${\mathbb Z} \times
{\mathbb Z}$, and let $\cal T$ be a triangulated $R$-category, and
let $M$ be an object in $\cal T$. Assume that there are two
invertible triangle functors $F$ and $G$ from $\cal T$ to itself
such that $FG$ is naturally isomorphic to $GF$ by $\delta: FG\ra
GF$, Suppose that $X\stackrel{\alpha}{\lra}
M_1\stackrel{\beta}{\lra} Y \stackrel{w}{\lra} X[1]$ is a triangle
in $\mathcal T$ such that $\alpha$ is a left $(\add(M), F, G,
\Phi)$-approximation of $X$ and $\beta$ is a right $(\add(M), F, G,
-(\Phi))$-approximation of $Y$. If $\Hom_{\cal T}(M,
G^iF^jX)=0=\Hom_{\cal T}(Y,G^iF^j(M))$ for $(0,0)\ne (i,j)\in \Phi$,
then $\Ex{\cal T}{F, G,\Phi}{X\oplus M}/I$ and $\Ex{\cal T}{F,
G,\Phi}{M\oplus Y}/J$ are derived equivalent. \label{with2functors}
\end{Theo}

Taking $G=id$, we recover Theorem \ref{verygthm1}. Taking
$G=[1]$ and $F=id$ yields a result on ``Ext-algebras". 
%To keep this
%note short, here we do not discuss any applications of Theorem
\ref{with2functors}.
\end{appendix}

{\footnotesize

\bigskip
Wei Hu, School of Mathematical Sciences, Laboratory of Mathematics
and Complex Systems, Beijing Normal University, 100875 Beijing,
People's Republic of  China

{\tt Email: huwei@bnu.edu.cn}

{\tt Homepage:  http://math.bnu.edu.cn/$\sim$huwei/}

\bigskip
Steffen Koenig, Institut f\"ur Algebra and Zahlentheorie,
Universit\"at Stuttgart, Pfaffenwaldring 57, 70569 Stuttgart,
Germany.

{\tt Email:  skoenig@mathematik.uni-stuttgart.de}

{\tt Homepage: http://www.iaz.uni-stuttgart.de/LstAGeoAlg/Koenig/}

\bigskip
Changchang Xi, School of Mathematical Sciences, Laboratory of
Mathematics and Complex Systems, Beijing Normal University, 100875
Beijing, People's Republic of  China

{\tt Email: xicc@bnu.edu.cn}

{\tt Homepage:  http://math.bnu.edu.cn/$\sim$ccxi/} }

%\bigskip March 23, 2010
\end{document}